\documentclass[a4paper,11pt]{elsarticle}
\usepackage{xcolor}
\usepackage{empheq}
\usepackage{float}
\usepackage{hyperref}
\usepackage{amsthm,amsmath,amscd,amsfonts,mathtools,stmaryrd,subcaption}
\usepackage[vmargin=2.5cm,hmargin=2.5cm,headheight=14.5pt,top=2cm,headsep=.5cm]{geometry}
\usepackage{mdframed}
\biboptions{sort&compress}

\newcommand{\irchi}{\mbox{\raisebox{1.5pt}{$\chi$}}}

\newtheoremstyle{ptheorem}{1em}{0em}{\itshape}{}{\bfseries}{.}{.5em}{\thmname{#1}\thmnumber{
		#2}\thmnote{ (\hspace{-.01pt}{#3})}}

\theoremstyle{ptheorem}

\newtheorem{thm}{Theorem}[section]
\newtheorem{pro}[thm]{Proposition}
\newtheorem{lem}[thm]{Lemma}
\newtheorem{cor}[thm]{Corollary}

\newtheoremstyle{hdef}{1em}{0em}{}{}{\bfseries}{.}{.5em}{\thmname{#1}\thmnumber{
		#2}\thmnote{ (\hspace{-.01pt}{#3})}}
\theoremstyle{hdef}

\newtheorem{dfn}[thm]{Definition}
\newtheorem{rem}[thm]{Remark}

\newtheorem{exa}[thm]{Example}

\numberwithin{equation}{section}
\numberwithin{figure}{section}
\allowdisplaybreaks

\journal{xx}

\begin{document}
	
	\begin{frontmatter}
		
		\title{The Wronskian and the variation of parameters method in the theory of linear Stieltjes differential equations of second order}
		
		\author{Francisco J. Fernández, Ignacio Marqu\'ez Alb\'es and F. Adri\'an F. Tojo
		}
		\address{Departamento de Estat\'{\i}stica, An\'alise Matem\'atica e Optimizaci\'on \\ Universidade de Santiago de Compostela \\ 15782, Facultade de Matem\'aticas, Campus Vida, Santiago, Spain. \\ e-mail: fjavier.fernandez@usc.es, ignacio.marquez@usc.es, fernandoadrian.fernandez@usc.es}
		
		\begin{abstract}
			
		In this work, we define the notions of Wronskian and simplified Wronskian for Stieltjes derivatives and study some of their properties in a similar manner to the context of time scales or the usual derivative. Later, we use these tools to investigate second order linear differential equations with Stieltjes derivatives  to find linearly independent solutions, as well as to derive the variation of parameters method  for problems with $g$-continuous coefficients. This theory is later illustrated with some examples such as the study of the one-dimensional linear Helmholtz equation with piecewise-constant coefficients.			
		\end{abstract}
		
		
	\end{frontmatter}
	
	
	\section{Introduction}
	
	The Wronskian has been a very useful tool in the theory of differential equations since its discovery \cite{HoeneWronski1812}. This instrument allows to check whether a given family of $n$ solutions of an order $n$ differential equation are linearly independent and, in fact, if the value of the Wronskian is known, it can also be used to find an $n$-th linearly independent solution of the differential equation of order $n$ provided $n-1$ linearly independent solutions are known.
	Wronskians are also important tools in the variation of parameters method, where the derivative of the parameters is expressed in terms of the Wronskian.
	
	These classical elements have also a role in newer theories of differential equations, such as the theory of Stieltjes differential equations, as this article will show. Unfortunately, the straightforward nature of the computations in the classical case (both regarding the Wronskian and the variation of parameters method) is not replicable in this setting due to the different nature of the product rule (see \cite[Proposition~3.9]{Fernandez2021}) which, in this case, reads 
	\[	\left(f_{1} f_{2}\right)_{g}'(t)=\left(f_{1}\right)_{g}'(t) f_{2}(t^*)+\left(f_{2}\right)_{g}'(t) f_{1}(t^*)+\left(f_{1}\right)_{g}'(t)\left(f_{2}\right)_{g}'(t) \Delta g(t^*),\]
	where $g$ is a nondecreasing and left-continuous function defining the Stieltjes derivative, $\Delta g$ denotes the jump of $g$ at a given point and $t^*$ is a point that depends on $t$. The authors have dedicated another work to explore in great detail the caveats and consequences of this product rule, see \cite{Fernandez2022b}.
	
	In this article we derive the definition of the Wronskian and the  variation of parameters method in the context of Stieltjes calculus, 
	taking into account the difficulties that arise from this theory and illustrating the applicability of the method with some examples. This endeavor will lead to the necessity of defining what we call a \emph{simplified Wronskian}, as well as the study of a special family of functions which, despite having low regularity, preserve the smoothness of a function when multiplied by it ---see Corollary~\ref{lemtec2}.
	
Given the existing relations between Stieltjes differential equations and other differential problems (see \cite[Section~8]{FriLo17}), our work draws from the classical theory as well as from the theory of differential equations in time scales, in which a notion of Wronskian also appears \cite{Bohner2001}.
	Remarkably, in this paper we will be able to apply our theory to the study of second order differential equations with non-constant coefficients, which cannot be found in \cite{Bohner2001}. Also, we would like to acknowledge that the notion of Wronskian for the case of Stieltjes differential equations also appears in the Master Thesis \cite{Lariviere}, although in that work it is not studied with an everywhere defined derivative, which is necessary for the finer points of the theory we develop below.
	
	To reach our goal, in Section 2 we give a brief overview of Stieltjes calculus in order to set the basis for our research, as well as introducing some new tools that are important for the work ahead, such as a version of the integration by parts formula or some results concerning the regularity of maps.
	In Section 3 we introduce the Wronskian and its simplified counterpart, presenting some of their basic properties. We apply the Wronskian and the variation of parameters method to the study second order Stieltjes differential equations in Section 4, and, finally, in Section 5, we provide three examples to which we apply the theory developed. In particular, we study the one-dimensional linear Helmholtz equation with piecewise-constant coefficients and we show that solution of the corresponding homogeneous equation arises naturally through the lenses of Stieltjes differential problems, which is particularly remarkable as this function is not two times continuously differentiable in the usual setting but it does present the corresponding regularity in the Stieltjes sense.
	
	\section{A brief overview of Stieltjes calculus}
	
	Let $g:\mathbb R\to\mathbb R$ be a nondecreasing and left-continuous function, which we shall refer to as a \emph{derivator}.
	We shall denote by 
	$\mu_g$ the Lebesgue-Stieltjes measure associated to $g$ given by
	\[ \mu_g([c,d))=g(d)-g(c),\quad c,d\in\mathbb R,\ c<d,\] 
	see \cite{Ru87,Sche97,Burk07}. We will use the term \emph{$g$-measurable} with respect to a set or function to refer to its $\mu_g$-measurability in the corresponding sense. Let $\mathcal L^1_{g}(X;\mathbb{C})$ the set of Lebesgue-Stieltjes $\mu_g$-integrable functions on a $g$-measurable set $X$ with values in $\mathbb{C}$, whose integral we denote by
\begin{displaymath}
	\int_X f(s)\,\operatorname{d}\mu_g(s),\quad f\in\mathcal L^1_{g}(X;\mathbb{C}).
\end{displaymath}

	Similarly, we will talk about properties holding \emph{$g$-almost everywhere} on a set $X$ (shortened to $g$-a.e. in~$X$), or holding for $g$-almost all (or simply, $g$-a.a.) $x\in X$, as a simplified way to express that they hold $\mu_g$-almost everywhere in $X$ or for $\mu_g$-almost all $x\in X$, respectively.

	Define the sets
	\begin{align*}
		C_g&=\{ t \in \mathbb R \, : \, \mbox{$g$ is constant on $(t-\varepsilon,t+\varepsilon)$ for some $\varepsilon>0$} \},\\
		D_g&=\{ t \in \mathbb R \, : \, \Delta g(t)>0\},
	\end{align*}
	where $\Delta g(t)=g(t^+)-g(t)$, $t\in\mathbb R$, and $g(t^+)$ denotes the right hand side limit of $g$ at $t$. First, observe that $C_g\cap D_g=\emptyset$. Furthermore, as pointed out in \cite{LoRo14}, the set $C_g$ is open in the usual topology of the real line, so it can be uniquely expressed as a countable union of open disjoint intervals, say
	\begin{equation}\label{Cgdisj}
		C_g=\bigcup_{n\in\Lambda} (a_n,b_n),
	\end{equation}
	where $\Lambda\subset \mathbb{N}$. With this notation, we introduce
	the sets $N_g^-$ and $N_g^+$ in \cite{LoMa19Resolution}, defined as
	\[ N_g^-=\{a_n: n\in\Lambda\}\setminus D_g,\quad N_g^+=\{b_n:n\in\Lambda\}\setminus D_g,\quad N_g=N_g^-\cup N_g^+.\] 
	
Before moving on to the study of the Stieltjes derivative, we present a result that contains some basic properties of the map $\Delta g$ that will be relevant for the work ahead. 

\begin{pro}[{\cite[Proposition~3.1]{Fernandez2022b}}]\label{regulated} For each $t\in\mathbb R$, we have that
\begin{equation}\label{Deltalimit}
\lim_{s\to t^-}\Delta g(t)=\lim_{s\to t^+}\Delta g(t)=0.
\end{equation}
In particular, $\Delta g$ is a regulated function, Borel-measurable and $g$-measurable.
\end{pro}

		After these considerations, we are finally in a position to define the Stieltjes derivative of a function on an interval as presented in \cite{Fernandez2021}, where the derivative is defined on the whole domain of the function, unlike in other papers such as \cite{FriLo17,LoRo14,Ma21}, where they exclude the points of $C_g$. This everywhere defined derivative is what will allow us to consider the second order derivative in a general setting, which is crucial, for instance, when $g$ is an step function, such as is the case when $g$ codifies the difference operator of difference equations. For a detailed discussion on the properties and implications of this everywhere defined derivative see \cite{Fernandez2021}.
		
		In order to present the derivative, we recall the hypotheses required in the mentioned paper, which we will also assume throughout this work. Namely, we will consider some $T>0$ and we will assume that $0\notin D_g\cup N_g^-$ and $T\notin N_g^+\cup D_g \cup C_g$. A careful reader might observe that throughout  \cite{Fernandez2021} it is also required that $g(0)=0$, however, this condition can easily be avoided by redefining the map $g$ if necessary, so we will not be considering it.

	\begin{dfn}[{\cite[Definition 3.7]{Fernandez2021}}]
		We define the \emph{Stieltjes derivative}\index{Stieltjes derivative}, or \emph{$g$-derivative}\index{$g$-derivative}, of a map $f:[a,b]\to\mathbb{C}$ at a point $t\in [0,T]$ as
		\[ 
		f'_g(t)=\left\{
		\begin{array}{ll}
			\displaystyle \lim_{s \to t}\frac{f(s)-f(t)}{g(s)-g(t)},\quad & t\not\in D_{g}\cup C_g,\vspace{0.1cm}\\
			\displaystyle\lim_{s\to t^+}\frac{f(s)-f(t)}{g(s)-g(t)},\quad & t\in D_{g},\vspace{0.1cm}\\
			\displaystyle\lim_{s\to b_n^+}\frac{f(s)-f(b_n)}{g(s)-g(b_n)},\quad & t\in C_{g},\ t\in(a_n,b_n),
		\end{array}
		\right.
		\] 
		where $a_n, b_n$ are as in~\eqref{Cgdisj}, provided the corresponding limits exist. In that case, we say that $f$ is \emph{$g$-differentiable at $t$}. 
	\end{dfn}

	\begin{rem}\label{remNgderivative}
		For $t\in N_g\cup\{0,T\}$,
		the corresponding limit in the definition of $g$-derivative at $t$ must be understood in the sense explained in \cite[Remark 2.2]{Ma21}, that is, the Stieltjes derivative at such points is computed as
		\begin{equation*}
			f'_g(t)=\left\{\begin{array}{ll}
				\displaystyle \displaystyle \lim_{s \to t^+}\frac{f(s)-f(t)}{g(s)-g(t)}, & t\in N_g^+\cup\{0\}, \vspace{0.1cm} \\
				\displaystyle \displaystyle \lim_{s \to t^-}\frac{f(s)-f(t)}{g(s)-g(t)}, & t\in N_g^-\cup\{T\},
			\end{array}\right.
		\end{equation*}
		provided the corresponding limit exists. Similarly, as pointed out in \cite{FriLo17}, $f$ is $g$-differentiable at $t\in D_g$ if and only if $f(t^+)$ exists and, in that case,
\begin{displaymath}		
f'_g(t)=\frac{f(t^+)-f(t)}{\Delta g(t)}.
\end{displaymath}
\end{rem}

\begin{rem}
		It is possible to further simplify the definition of the Stieltjes derivative at a point $t\in[a,b]$ by defining 
		\begin{equation*}
			t^*=
			\begin{dcases}
				t,\quad & t\not\in C_g,\vspace{0.1cm}\\
				b_n,\quad & t\in (a_n,b_n)\subset C_g,
			\end{dcases}
		\end{equation*}
		with $a_n,b_n$ as in~\eqref{Cgdisj}. With this notation, we have that 
		\[ 
		f'_g(t)=\left\{
		\begin{array}{ll}
			\displaystyle \lim_{s \to t}\frac{f(s)-f(t)}{g(s)-g(t)},\quad & t\not\in D_{g}\cup C_g,\vspace{0.1cm}\\
			\displaystyle\lim_{s\to t^{*+}}\frac{f(s)-f(t^*)}{g(s)-g(t^*)},\quad & t\in D_{g}\cup C_g,
		\end{array}
		\right.
		\] 
		provided the corresponding limit exists. Note that the information in Remark~\ref{remNgderivative} should still be taken into account. From now on, given a function $f$, we denote by $f^*$ the function defined as $f^*(t)=f(t^*)$. For instance, $\Delta g^*(t)=\Delta g(t^*)$. 
		Observe that $\{t\in [0,T]:\; t\neq t^*\}\subset C_g$, therefore, 
		\begin{equation}\label{medtneqts}
		0\leq\mu_g(\{t\in [0,T]:\; t\neq t^*\})\leq \mu_g(C_g)=0.
		\end{equation}
	\end{rem}

	The following result,  \cite[Proposition~3.9]{Fernandez2021},  includes some basic properties of this derivative.

	\begin{pro}\label{PropStiDer} 
		Let $t\in[a,b]$.	 If $f_1,f_2:[0,T]\to\mathbb{C}$ are $g$-differentiable at $t$, then:
		\begin{itemize}
			\item The function $\lambda_{1} f_{1}+\lambda_{2} f_{2}$ is $g$-differentiable at $t$ for any $\lambda_{1}, \lambda_{2} \in \mathbb{C}$ and
			\begin{equation*}%
				\left(\lambda_{1} f_{1}+\lambda_{2} f_{2}\right)_{g}'(t)=\lambda_{1}\left(f_{1}\right)_{g}'(t)+\lambda_{2}\left(f_{2}\right)_{g}'(t).
			\end{equation*}
			\item\emph{ (Product rule).} The product $f_{1} f_{2}$ is $g$-differentiable at $t$ and
			\begin{equation*}
				\left(f_{1} f_{2}\right)_{g}'(t)=\left(f_{1}\right)_{g}'(t) f_{2}(t^*)+\left(f_{2}\right)_{g}'(t) f_{1}(t^*)+\left(f_{1}\right)_{g}'(t)\left(f_{2}\right)_{g}'(t) \Delta g(t^*).
			\end{equation*}
			\item \emph{ (Quotient rule).} If $f_2(t^*)\,(f_2(t^*)+(f_2)'_g(t)\, \Delta g(t^*))\neq 0$, the quotient $f_1/f_2$ is $g$-differentiable at $t$ and 
			\begin{equation*}
				\left(\frac{f_1}{f_2}\right)'_g(t)=\frac{\left(f_{1}\right)_{g}'(t) f_{2}(t^*)-\left(f_{2}\right)_{g}'(t) f_{1}(t^*)}{f_2(t^*)\,(f_2(t^*)+(f_2)'_g(t)\, \Delta g(t^*))}.
			\end{equation*}
		\end{itemize}
	\end{pro}

The following result, presents some conditions ensuring that the map $\Delta g^*$ in the product rule is $g$-differentiable. This will be  a fundamental result for the variation of parameters method.  
\begin{pro}[{\cite[Corollary 4.4]{Fernandez2022b}}]\label{diffDeltag*}
		Consider the sets 
		\begin{align*}
			D_1&=\{t\in[0,T]\cap N_g^-: t\in (D_g\cap[0,t))'\},\\
			D_2&=\{t\in[0,T]\cap N_g^+: t\in (D_g\cap(t,T])'\},\\
			D_3&=\{t\in[0,T]\backslash (N_g\cup D_g: t\in (D_g\cap[0,T])'\},
		\end{align*}
		and assume  that
		\begin{equation}\label{condDeltaderivable}
				\lim_{\substack{s\to t\\ s\in D_g}}\frac{\Delta g|_{[0,T]}(s)}{g(s)-g(t)}=0,\quad\mbox{for all }
				t\in D_1\cup D_2\cup D_3.
		\end{equation}
		Then, $\Delta g^*|_{[0,T]}$ is $g$-differentiable on $[0,T]$ and $\left(\Delta g^*\right)'_g=-\irchi_{D_g}^*$, where $\irchi_{D_g}$ denotes the indicator function on $D_g$.
\end{pro}	

		In the context of Stieltjes calculus, we also find a concept of continuity related to the map $g$.
		\begin{dfn}[{\cite[Definition 3.1]{FriLo17}}]\label{dfncont} A function $f:[0,T]\to {\mathbb C}$ is
			\emph{$g$-continuous} at a point $t\in [0,T]$,
			or \emph{continuous with respect to $g$} at $t$, if for every $\varepsilon>0$, there exists $\delta>0$ such that 
			\[|f(t)-f(s)|<\varepsilon,\quad \mbox{for all }s\in[0,T],\ |g(t)-g(s)|<\delta.\]
			If $f$ is $g$-continuous at every point $t\in [0,T]$, we say that $f$ is $g$-continuous on $[0,T]$. We denote by $\mathcal{C}_g([0,T];\mathbb{C})$ the \emph{set of $g$-continuous functions} on $[0,T]$; and  by $\mathcal{BC}_g([0,T];\mathbb{C})$ 
			the \emph{set of bounded $g$-continuous functions} on $[0,T]$
		\end{dfn}
		\begin{rem}
			Observe that we distinguish between $\mathcal{C}_g([0,T];\mathbb{C})$ and $\mathcal{BC}_g([0,T];\mathbb{C})$ because, as pointed out in \cite[Example 3.19]{MarquezTesis}, $g$-continuous function on compact intervals need not be bounded. It is important to note that, as explained in \cite[Example 3.23]{MarquezTesis}, $g$-differentiable functions need not be $g$-continuous either.
		\end{rem}

		\begin{rem}
			The set $\mathcal{BC}_g([0,T];\mathbb{C})$ equipped with the supremum norm, $\left\lVert \cdot\right\rVert_0$ is a Banach space. As such, it is possible to talk about linearly indepence in $\mathcal{BC}_g([0,T];\mathbb{C})$. We shall say that $v_1,v_2 \in \mathcal{BC}_g([0,T];\mathbb{C})$ are \emph{linearly dependent} if there exist $c_1,c_2\in\mathbb C$ such that
			\[c_1 v_1(t)+c_2 v_2(t)=0,\quad  t\in [0,T].\] 
			Otherwise, we say that $v_1$ and $v_2$ are \emph{linearly independent}.
		\end{rem}
		
Naturally, we have that the sum and product of $g$-continuous functions are $g$-continuous. Similarly, the quotient of two $g$-continuous functions is also $g$-continuous provided that the function on the denominator does not vanish, see \cite[Lemma 2.14.]{Fernandez2022b}.
		The following result describes some other basic properties for $g$-continuous functions. It can be  deduced directly from \cite[Proposition~3.2]{FriLo17}, in which the same information is presented for real-valued functions.
		
		\begin{pro}\label{proreg} 
			If $f:[0,T] \to \mathbb{C}$ is $g$-continuous on $[0,T],$ then:
			\begin{itemize}
				\item $f$ is continuous from the left at every $t \in(0,T]$;
				\item if $g$ is continuous at $t \in[0,T),$ then so is $f$;
				\item if $g$ is constant on some $[\alpha, \beta] \subset[0,T],$ then so is $f$.
			\end{itemize}
		\end{pro}

\begin{rem}\label{f*equalsf} As a direct consequence of Proposition~\ref{proreg}, we see that if $f\in\mathcal{C}_g([0,T];\mathbb{C})$, then $f^*=f$. Indeed, let $t\in [0,T]$ and let us show that $f(t)=f^*(t)$. If $t\not\in C_g$ this is trivial as $t^*=t$. For $t\in C_g$, $g$ is constant on $[t,t^*]$ and, thus, so is $f$ so, $f(t)=f(t^*)=f^*(t)$.
			
In particular, this means that if $f_1,f_2\in\mathcal{C}_g([0,T];\mathbb{C})$ are $g$-differentiable at $t$, then 
\begin{displaymath}
\left(f_{1} f_{2}\right)_{g}'(t)=\left(f_{1}\right)_{g}'(t) f_{2}(t)+\left(f_{2}\right)_{g}'(t) f_{1}(t)+\left(f_{1}\right)_{g}'(t)\left(f_{2}\right)_{g}'(t) \Delta g(t^*).
\end{displaymath}
\end{rem}	
		
Next, we introduce the concept of $g$-absolute continuity, which is the extension of the notion of absolute continuity to the context of Stieltjes calculus. This definition was presented in \cite{LoRo14} and it connects the Stieltjes derivative and the Lebesgue-Stieltjes integral, see  \cite[Proposition~5.4]{LoRo14}. Here, we introduce its definition as part of the mentioned result, which was originally stated for real-valued functions but easily extends to complex-valued ones.
		
		\begin{thm}\label{gFTC2}
			Let $F:[0,T]\to\mathbb C$. The following conditions are equivalent:
			\begin{enumerate}
				\item[\textup{1.}] The function $F$ is \emph{$g$-absolutely continuous on $[0,T]$}, according to the following definition: for every $\varepsilon>0$, there exists $\delta>0$ such that for every open pairwise disjoint family of subintervals $\{(a_n,b_n)\}_{n=1}^m$,
				\begin{displaymath}
					\sum_{n=1}^m (g(b_n)-g(a_n))<\delta
					\implies
					\sum_{n=1}^m |F(b_n)-F(a_n)|<\varepsilon.
				\end{displaymath}
				\item[\textup{2.}] The function $F$ satisfies the following conditions:
				\begin{enumerate}
					\item[\textup{(i)}] there exists $F'_g(t)$ for $g$-a.a. $t\in[0,T)$;
					\item[\textup{(ii)}] $F'_g\in \mathcal L^1_{g}([0,T),\mathbb C)$;  
					\item[\textup{(iii)}] for each $t\in[0,T]$, 
					\begin{equation}\label{Barrow}
						F(t)=F(0)+\int_{[0,t)}F'_g(s)\,\operatorname{d}\mu_g(s).
					\end{equation}
				\end{enumerate}
			\end{enumerate}
		\end{thm}
		We denote by  $\mathcal{AC}_g([0,T];\mathbb{C})$ the \emph{set of $g$-absolutely continuous functions} on $[0,T]$. 
		
		The following result is a version of the formula of integration by parts for $g$-absolutely continuous functions, which is a direct consequence of  Proposition~\ref{PropStiDer} and Theorem~\ref{gFTC2}.
		
		\begin{lem}[Integration by parts]\label{InByPart} Given $w_1,w_2\in 
			\mathcal{AC}_g([0,T];\mathbb{C})$, we have that 
			$w_1w_2 \in \mathcal{AC}_g([0,T];\mathbb{C})$ and, furthermore, for each $t\in[0,T]$,
			\begin{displaymath}
				\begin{aligned}
					w_1(t)\, w_2(t)-w_1(0)\, w_2(0) =& \int_{[0,t)} (w_1)'_g\, w_2\, \operatorname{d}\mu_g+
					\int_{[0,t)} w_1\, (w_2)'_g\, \operatorname{d}\mu_g+
					\int_{[0,t)} (w_1)'_g\, (w_2)'_g\, \Delta g\, \operatorname{d}\mu_g.
				\end{aligned}
			\end{displaymath}
		\end{lem}
		
		\begin{proof} First, observe that \cite[Proposition 5.4]{FriLo17} 
			ensures that $w_1w_2 \in \mathcal{AC}_g([0,T];\mathbb{C})$. 
			Now, Theorem~\ref{gFTC2} and Remark~\ref{f*equalsf} ensure that
			\begin{displaymath}
				(w_1 w_2)'_g=w_1\, (w_2)'_g+ (w_1)'_g\, w_2 + 
				(w_1)'_g \, (w_2)'_g\, \Delta g^*,\quad g\mbox{-a.e. in }[0,T),
			\end{displaymath}
			so, thanks to~\eqref{medtneqts}, we obtain
			\begin{displaymath}
				(w_1 w_2)'_g=w_1\, (w_2)'_g+ (w_1)'_g\, w_2 + 
				(w_1)'_g \, (w_2)'_g\, \Delta g,\quad g\mbox{-a.e. in }[0,T),
			\end{displaymath}
			from which the integration by parts formula follows using \eqref{Barrow}.
		\end{proof}
		
		As pointed out in  \cite[Proposition 5.5]{FriLo17}, $\mathcal{AC}_g([0,T];\mathbb{C})
		\subset \mathcal{BC}_g([0,T];\mathbb{C})$ so every $g$-absolutely continuous function is $g$-continuous and, as such, presents the properties introduced before. Note, however, that $g$-absolutely continuous functions are not, in general, $g$-differentiable everywhere, which motives the following definition, see \cite[Definitions~3.11 and 3.12]{Fernandez2021}.
		
		\begin{dfn} 
			Given $k\in\mathbb N$, we define $\mathcal{C}^0_g([0,T];{\mathbb C}):=\mathcal{C}_g([0,T];{\mathbb C})$ and $\mathcal{C}^k_g([0,T];{\mathbb C})$ recursively as
			\begin{displaymath}
				\mathcal{C}^k_g([0,T]):=\{f \in \mathcal{C}^{k-1}([0,T];{\mathbb C}):\; (f_g^{(k-1)})'_g\in
				\mathcal{C}_g([0,T];{\mathbb C})\},\end{displaymath}
			where $f^{(0)}_g:=f$ and $f^{(k)}_g:=(f^{(k-1)}_g)'_g$, $k\in\mathbb N$.
			Similarly, given $k\in\mathbb N$, we define $\mathcal{BC}^0_g([0,T];{\mathbb C}):=\mathcal{BC}_g([0,T];{\mathbb C})$ and $\mathcal{BC}^k_g([0,T];{\mathbb C})$ recursively as
			\begin{displaymath}
				\mathcal{BC}_g^k([a,b];{\mathbb C}):=\{f \in \mathcal{C}_g^k([a,b];{\mathbb C}):\; f^{(n)}_g\in \mathcal{BC}_g([a,b];{\mathbb C}),\;
				\forall n=0,\ldots,k\}.
			\end{displaymath}
			We also define $\mathcal{C}^\infty_g([0,T];{\mathbb C}):=\bigcap_{{n\in{\mathbb N}}}\mathcal{C}^k_g([0,T];{\mathbb C})$ and
			$\mathcal{BC}^\infty_g([0,T];{\mathbb C}):=\bigcap_{{n\in{\mathbb N}}}\mathcal{BC}^k_g([0,T];{\mathbb C})$. 
		\end{dfn}

One of the most important examples of functions in the space 
$\mathcal{BC}_g^{\infty}([a,b];\mathbb{C})$ is the $g$-exponential of a constant function. In the following definition, we collect some of the information available on \cite{FriLo17, Ma21, Fernandez2021} for the $g$-exponential of an integrable function.
\begin{dfn}\label{gexpdef}
	Given a function $p:[0,T]\rightarrow \mathbb{C}$, we say that it is \emph{regressive} if
	\begin{equation}\label{pcond}
		1+p(t)\,\Delta g(t)\not=0,\quad t\in [0,T]\cap D_g.
	\end{equation}

	Given a regressive function $p\in \mathcal{L}_g^1([0,T);\mathbb{C})$, we define the \emph{$g$-exponential} associated to the map $p$ as 
	\begin{equation}\label{eq:gexpotential}
		\exp_g(p;t):=
		\exp \left( \int_{[0,t)} \tilde{p}(s)\, \operatorname{d}\mu_g(s)\right),\quad t\in[0,T],
	\end{equation}
	where, denoting by $\operatorname{ln}(z):=\operatorname{ln}\left\lvert z\right\rvert +i\operatorname{Arg}(z)$ the principal 
	branch of the complex logarithm, 
	\begin{displaymath}
		\tilde{p}(s):=\begin{dcases}
			p(s), & s \in [0,T)\setminus D_g, \\
			\frac{\operatorname{ln}\left(1+p(s)\, \Delta g(s)\right)}{\Delta g(s)}, & s \in [0,T)\cap D_g.
		\end{dcases}
	\end{displaymath}
\end{dfn}
\begin{rem}\label{expCinfty}
	The $g$-exponential map belongs to $\mathcal{AC}_g([0,T];\mathbb{C})$, see \cite[Theorem 4.2]{Fernandez2021}, and, furthermore, it is the only function in that space satisfying
	\begin{equation} \label{eqgexp}
		\left\{\begin{aligned}
			&v'_g(t)=p(t)\,  v(t),\; g\mbox{-a.a.}\, t \in [0,T], \\
			&v(0)=1.
		\end{aligned}\right.
	\end{equation}
In particular, if $p\in 
\mathcal{BC}_g([0,T];\mathbb{C})$ then $\exp_g(p;\cdot)\in\mathcal{BC}^1_g([0,T];\mathbb{C})$. Furthermore, if $p(t)=\lambda\in\mathbb C$, $t\in[0,T]$, then $\exp_g(p;\cdot)\in\mathcal{BC}^\infty_g([0,T];\mathbb{C})$.
\end{rem}
	
	In order to make this work more self-contained, we highlight below some of the properties of the $g$-exponential function whose proof can be found in \cite[Proposition~4.6]{Fernandez2021}.
	
\begin{pro}\label{gexpprop} Let $p,q\in \mathcal{L}_g^1([0,T);\mathbb{C})$ be two regressive functions. Then: 
\begin{enumerate}
\item For each $t\in[0,T]$, 
\begin{displaymath}
\exp_g(p;t)\, \exp_g(q;t)=
\exp_g(p+q+p\, q \, \Delta g;t).
\end{displaymath}
In particular,
\begin{displaymath}
\exp_g(p;t)^2=\exp_g(2\,p+p^2\, \Delta g;t),\quad t\in [0,T].
\end{displaymath}
\item For each $t\in [0,T]$, 
\begin{displaymath}
\exp_g(p;t)^{-1}=\exp_g\left( 
-\frac{p}{1+p\, \Delta g};t\right).
\end{displaymath}
As a consequence,
\begin{displaymath}
\exp_g(p;t)\, \exp_g(q;t)^{-1}=
\exp_g\left(\frac{p-q}{1+q\, \Delta g};t \right),\quad t \in [0,T].
\end{displaymath}
\end{enumerate} 
\end{pro}
	
One of the main problems that we will encounter when trying to apply the method of variation of parameters in the context of Stieltjes calculus is the fact that the product of two functions in the space 
$\mathcal{BC}_g^1([0,T];\mathbb{C})$ 
is not necessarily a function in the space 
$\mathcal{BC}_g^1([0,T];\mathbb{C})$ as a consequence of the expression of the product rule given by Proposition~\ref{PropStiDer}
(see \cite[Remark~3.16]{Fernandez2021} and, for a detailed 
discussion, \cite{Fernandez2022b}). Nevertheless, for a given function $f\in \mathcal{BC}_g^1([0,T];\mathbb{C})$, it is possible, in some cases, to find another function, $h\in \mathcal{BC}_g([0,T];\mathbb{C})$, $g$-differentiable on $[0,T]$, such that the product $fh$ lays in $\mathcal{BC}_g^1([0,T];\mathbb{C})$. Indeed, under these assumptions, Remark~\ref{f*equalsf} ensures that
	\[(fh)'_g=f_g'\, h+f\, h_g'+f_g'\, h_g'\, \Delta g^*=f_g'\, h+h_g'\left[f+f_g'\, \Delta g^* \right].\]
	Therefore, if $f+f_g'\, \Delta g^*\not=0$ on $[0,T]$, any function $h$ such that 
	\begin{equation}\label{eq:regderpro}
		h_g'=\frac{\eta}{f+f_g'\, \Delta g^*},
	\end{equation}
	with $\eta \in \mathcal{BC}_g([0,T];\mathbb{C})$,  would yield that $(fh)'_g=f_g'\, h+\eta$, which would belong to $\mathcal{BC}_g^1([0,T];\mathbb{C})$, as we wanted. Note that, in general, we still would not have that $h\in\mathcal{BC}_g^1([0,T];\mathbb{C})$, as $h_g'$ might  not belong to $\mathcal{BC}_g([0,T];\mathbb{C})$. This reasoning is the idea behind Corollary~\ref{lemtec2}, which we present as a consequence of the following technical result.

\begin{lem} \label{dercont} Let $f:[0,T]\rightarrow \mathbb{C}$ be a bounded function which is continuous (in the usual sense) at every $t \in [0,T]\setminus D_g$. Then the map $\psi:[0,T]\to\mathbb C$ defined as 
\begin{displaymath}
\psi(t)=\int_{[0,t)} f(s)\, \operatorname{d}\mu_g(s),\quad t\in[0,T],
\end{displaymath}
is well-defined;
belongs to $\mathcal{AC}_g([0,T];\mathbb{C})$ and
\begin{equation}\label{eq:lemtec1bis}
\psi'_g(t)=f^*(t),\quad \mbox{for all }t\in[0,T].
\end{equation}
\end{lem}

\begin{proof} 
First, observe that $f\in \mathcal{L}_g^1([0,T);\mathbb{C})$. 
Indeed, since $f$ is continuous 
on $[0,T]\setminus D_g$ and $D_g$ is a countable set, 
we have that $f$ is Borel measurable, thus $g$-measurable. Now, the $g$-integrability is clear since $f$ is bounded. Thus, Theorem~\ref{gFTC2} 
ensures that
$\psi \in \mathcal{AC}_g([0,T];\mathbb{C})$ and 
$\psi_g'(t)=f(t)$ for  $g$-a.a. 
$t\in [0,T]$ and, in particular, for every $t\in [0,T)\cap D_g$. 
Hence, we need to show that equation \eqref{eq:lemtec1bis} holds on
$[0,T]\setminus D_g$.

First, let $t\notin C_g \cup D_g$. In this case, we reason similarly to \cite[Lemma 3.14]{Fernandez2021}; namely, computing the limit 
\begin{displaymath}
	\lim_{s\rightarrow t} \frac{\psi(s)-\psi(t)}{g(s)-g(t)}
\end{displaymath}
on the domain of the function, that is, $A_t:=\{s\in [0,T]:\; g(s)\neq g(t)\}$. 
Let $\varepsilon>0$. Since $f$ is continuous on $[0,T]\setminus D_g$, there exits $\delta>0$ such that $\left\lvert f(u)-f(t)\right\rvert <\varepsilon$ if $\left\lvert u-t\right\rvert <\delta$. Now, 
for $s\in [0,T]\cap A_t$ such that $\left\lvert t-s\right\rvert <\delta$, denoting $[\left\lvert t,s\right\rvert ):=[\min\{t,s\},\max\{t,s\})$, we have that
\begin{displaymath}
	\begin{aligned}
		\left|\frac{\psi(s)-\psi(t)}{g(s)-g(t)}-f(t)
		\right| =& \left|
		\frac{\operatorname{sgn}(s-t)}{g(s)-g(t)}
		\int_{[|t,s|)} f(u)\, \operatorname{d}\mu_g(u)-f(t)
		\right| \\
		\leq &\frac{1}{|g(s)-g(t)|} \int_{[|t,s|)}
		\left| f(u)-f(t)\right|\, \operatorname{d}\mu_g(u) \leq \varepsilon.
	\end{aligned}
\end{displaymath}
Thus, since $t\notin C_g \cup D_g$, it follows that $t^*=t$, so
\begin{displaymath}
	\lim_{s\rightarrow t} \frac{\psi(s)-\psi(t)}{g(s)-g(t)}=f(t)=f(t^*),
\end{displaymath}
as we wanted.

Finally, if $t\in C_g$, then $t^*\in D_g\cup N_g^+$, and we already know that $\psi$ is $g$-differentiable at $t^*$ and $\psi'_g(t^*)=f^*(t^*)=f^*(t)$. Hence, by the definition of the $g$-derivative at a point of in $C_g$, we have that $\psi'_g(t)=\psi'_g(t^*)=f(t^*)$, which finishes the proof of the result.
\end{proof}

\begin{rem} Observe that in Lemma~\ref{dercont}, even though we can only ensure that $\psi\in\mathcal{AC}_g([0,T];\mathbb C)$, we can still guarantee that the function $\psi$ is $g$-differentiable on the whole $[0,T]$.
\end{rem}

As anticipated, we have the following corollary which, for a given function, provides an explicit expression for another function such that their product is an element of $\mathcal{BC}_g^1([0,T];\mathbb{C})$.

\begin{cor} \label{lemtec2} Let $\eta\in \mathcal{BC}_g([0,T];\mathbb{C})$ and 
$f\in \mathcal{BC}_g^1([0,T];\mathbb{C})$ be such that the function
\begin{displaymath}
\widetilde{f}(t):=\frac{1}{f(t)+
f_g'(t)\, \Delta g(t)},\quad t\in[0,T],
\end{displaymath}
is well-defined and bounded. Then, the map $\varphi:[0,T]\to\mathbb C$ defined as 
\begin{displaymath}
\varphi(t):=\int_{[0,t)} \widetilde{f}(s)\,\eta(s)\,  \operatorname{d}\mu_g(s),\quad t\in[0,T],
\end{displaymath}
is well-defined; belongs to $\mathcal{AC}_g([0,T];\mathbb{C})$ and $\varphi'_g(t)= \widetilde f(t^*)\,\eta(t)$, $t\in[0,T]$. Furthermore,
\begin{displaymath}
(f\varphi)'_g(t)=f_g'(t) \varphi(t)+ \eta(t),\quad t\in[0,T],
\end{displaymath}
and, as a consequence, $f\varphi\in \mathcal{BC}_g^1([0,T];\mathbb{C})$.
\end{cor}

\begin{proof} 
	Observe that it is enough to show that $\widetilde{f}\eta$ satisfies the conditions of Lemma~\ref{dercont}, namely, it is bounded and continuous at every $t\in[0,T]\setminus D_g$. Note that the boundedness follows directly from the hypotheses.

	Let $t\in [0,T]\setminus D_g$ and $\{t_n\}_{n \in \mathbb{N}}\subset [0,T]$ be a sequence such that $t_n \rightarrow t$. In that case, Proposition~\ref{regulated} ensures that $\Delta g(t_n)\rightarrow 0$, so
\[
\lim_{n \to \infty} \widetilde{f}(t_n)= \lim_{n\to \infty} 
\frac{1}{f(t_n)+f_g'(t_n)\, \Delta g(t_n)} 
= \frac{1}{f(t)}=\widetilde{f}(t),
\]
where we have used the fact that $f$ is continuous at $t$ since $f\in \mathcal{BC}_g^1([0,T];\mathbb{C})$ and $t\not\in D_g$, see 
Proposition \ref{proreg}. This shows that $\widetilde f$ is continuous at $t$. Now Proposition \ref{proreg} guarantees, once again that $\eta$ is continuous at $t$ so $\widetilde f\eta$ is continuous at $t$.

Hence, we can apply Lemma~\ref{dercont} to see that $\varphi$ is well-defined, belongs to $\mathcal{AC}_g([0,T];\mathbb{C})$ and 
\begin{displaymath}
\varphi'_g(t)= \widetilde f(t^*)\,\eta(t)=
\frac{\eta(t)}{f(t)+f_g'(t)\,\Delta g(t^*)},\; t\in[0,T].
\end{displaymath}
Finally, Proposition~\ref{PropStiDer} ensures that, for each $t\in[0,T]$
\begin{displaymath}
\begin{aligned}
	(f\, \varphi)'_g(t)=& f_g'(t)\, \varphi(t)+f(t)\, \varphi'_g(t)+
	f_g'(t)\, \varphi'_g(t)\, \Delta g(t^*)\\
	=& f'_g(t)\, \varphi(t)+ \varphi_g'(t)\left(f(t)+f_g'(t)\, \Delta g(t^*)\right)= f_g'(t)\,\varphi(t)+\eta(t),
\end{aligned}
\end{displaymath}
which finishes the proof.
\end{proof}

\begin{rem}\label{remboundedfromzero}
	In Corollary~\ref{lemtec2}, we are assuming that $\widetilde f$ is well-defined and bounded, which is a technical condition required for all the functions to be well-defined and some of the conditions of Lemma~\ref{lemtec2} to be satisfied. However, if $f$ is bounded away from zero, those conditions can be dropped, given a condition which might be easier to check. Indeed, suppose $f$ is bounded away from zero. Then, there exists $M>0$ such that $|f(t)|>M,\quad t\in[0,T]$. 
	Now, since
	$f\in \mathcal{BC}_g^1([0,T];\mathbb{C})$, we have as a direct consequence of Remark~\ref{remNgderivative} and Proposition~\ref{proreg} that $f(t)+f_g'(t)\, \Delta g(t)=f(t^+)$, $t\in [0,T]$,
	so, 
	\begin{equation*}\label{ftildewelldefined}
		|f(t)+f_g'(t)\, \Delta g(t)|=|f(t+)|\ge M>0,\quad t\in[0,T],
	\end{equation*}
	which guarantees that $\widetilde f$ is well-defined and bounded on $[0,T]$.
\end{rem}

\begin{cor} \label{lemtec1} Let $\eta \in \mathcal{BC}_g([0,T];\mathbb{C})$ 
and $\lambda \in \mathbb{C}$ be 
such that $1+\lambda \Delta g(t)\neq 0$, 
$t\in [0,T]\cap D_g$. Then, the map $\varphi:[0,T]\to\mathbb C$ defined as
\begin{displaymath}
\varphi(t): = \int_{[0,t)} \frac{\eta(s)}{1+\lambda \Delta g(s)}\, \operatorname{d}\mu_g(s),\quad t\in[0,T],
\end{displaymath}
is well-defined; belongs to $\mathcal{AC}_g([0,T];\mathbb{C})$ and
\begin{displaymath}
\varphi'_g(t)=\frac{\eta(t)}{1+\lambda \Delta g^*(t)},\quad t\in[0,T].
\end{displaymath}
\end{cor}
\begin{proof} Observe that, in order to obtain the result, it is enough to show that the map 
\begin{displaymath}
\theta(t)=\frac{\eta(t)}{1+\lambda \Delta g(t)},\quad t\in [0,T],
\end{displaymath}
satisfies the hypotheses of Lemma~\ref{dercont}.

We start by proving that $\theta$ is bounded. To this end, define
\begin{equation*}
A=\{t\in [0,T):\; \left\lvert 1+\lambda \Delta g(t) \right\rvert \geq 1/2\}, \quad
B=\{t \in [0,T):\; \left\lvert 1+\lambda \Delta g(t)\right\rvert <1/2\}.
\end{equation*}
Observe that, necessarily, $B\subset[0,T)\cap D_g$. We claim that $B$ has finite 
cardinality. Indeed, since 
\begin{displaymath}
0\le\sum_{t\in [0,T)\cap D_g} \Delta g(t)=\int_{[0,T)\cap D_g}\, \operatorname{d}\mu_g(u)\le \int_{[0,T)}\, \operatorname{d}\mu_g(u)=\mu_g([0,T))=g(T)-g(0)<\infty,
\end{displaymath}
the set $X:=\{[0,T)\cap D_g :\; |\lambda| \Delta g(t)>1/2\}$ is finite. Now, given an element $t\in B$, we have that if $1-|\lambda|
\Delta g(t)<|1+\lambda \Delta g(t)|<1/2$, then $t\in X$, so $B$ 
is finite. Hence,
\[
\left\lVert \theta\right\rVert_0\leq 
\sup_{t\in A}\left\lVert \theta(t)\right\rVert + 
\sup_{t\in B}\left\lVert \theta(t)\right\rVert\leq 
\left\lVert  \eta\right\rVert_0 \left( 
2+\max\left\{\left\lvert 1+\lambda \Delta g(t)\right\rvert ^{-1}:\, 
t\in B\right\}
\right)<\infty.
\]

Now, the proof of the continuity of $\theta$ on $[0,T]\setminus D_g$ is analogous to the proof of the continuity of $\widetilde{f}$ in the Corollary~\ref{lemtec2}, so we omit it.

Therefore, we have that the hypotheses of Lemma~\ref{dercont} are satisfied, so 
$\varphi\in \mathcal{AC}_g([0,T];\mathbb{C})$ 
and 
\begin{displaymath}
\varphi'_g(t)=\frac{\eta(t^*)}{1+\lambda \Delta g(t^*)}=
\frac{\eta(t)}{1+\lambda \Delta g(t^*)},\quad t\in[0,T],
\end{displaymath}
where the last equality is a consequence of the $g$-continuity of $\eta$, see Remark~\ref{f*equalsf}.
\end{proof}

Finally, we include a result that will be of interest for the work ahead and it can be directly obtained from Corollary~\ref{lemtec1}. 

\begin{cor} \label{prodcinf} 
	Let $\lambda \in \mathbb{C}$ be 
such that $1+\lambda \Delta g(t)\neq 0$, $t\in [0,T]\cap D_g$, and define
\begin{displaymath}
\varphi(t)=\int_{[0,t)} \frac{1}{1+\lambda \Delta g(s)}\, \operatorname{d}\mu_g(s),\;
t\in [0,T].
\end{displaymath}
Then, the map $v(t)=\varphi(t)\, \exp_g(\lambda;t)$,  $t\in[0,T]$, belongs to $\mathcal{BC}_g^{\infty}([0,T];\mathbb{C})$ and, moreover,
\begin{equation} \label{eq:dern}
v_g^{(n)}(t)= n\, \lambda^{n-1} \exp_g(\lambda;t)+ \lambda^n\, v(t),\quad t \in [0,T],\quad n\in\mathbb N.
\end{equation}

\end{cor}

\begin{proof} 
	First, observe that Corollary~\ref{lemtec1} and \cite[Theorem~4.2]{Fernandez2021}  are enough to guarantee that $v\in \mathcal{BC}_g([0,T];\mathbb{C})$. Hence, it suffices to prove that \eqref{eq:dern} holds as, in that case, it follows that $v_g^{(n)}\in\mathcal{BC}_g([0,T];\mathbb{C})$ as it can be expressed in terms of $v$ and the $g$-exponential. We prove \eqref{eq:dern}  by induction on $n\in\mathbb N$.

Let $t\in[0,T]$. Since $\varphi$ and $\exp_g(\lambda;\cdot)$ are $g$-differentiable at $t$,
Proposition~\ref{PropStiDer}, Corollary~\ref{lemtec1}  and Remark~\ref{f*equalsf} yield
\begin{displaymath}
\begin{aligned}
v'_g(t)=&\lambda \, \varphi(t^*)\exp_g(\lambda;t)+\varphi_g'(t) \exp_g(\lambda;t^*)+
\lambda \varphi_g'(t) \exp_g(\lambda;t)\Delta g(t^*)\\
=& \lambda\, \varphi(t) \exp_g(\lambda;t) + \exp_g(\lambda;t),
\end{aligned}
\end{displaymath}
which proves the case $n=1$. Assume now that \eqref{eq:dern} holds for every $n\in\mathbb N$ and let us show that it holds for $n+1$. Observe that this means that $v_g^{(n)}\in\mathcal{BC}_g([0,T];\mathbb{C})$. Let $t\in[0,T]$. In that case, since $v_g^{(n)}$ and $v$ are $g$-differentiable at $t$, Proposition~\ref{PropStiDer} and Remark~\ref{f*equalsf} ensure that
\begin{displaymath}
\begin{aligned}
v_g^{(n+1)}(t)&=v_g^{(n)}(t)v(t)+v_g^{(n-1)}(t)v'_g(t)+v_g^{(n)}(t)v'_g(t)\Delta g(t^*)\\
 &=k\, \lambda^{n} \exp_g(\lambda;t)+ \lambda^n\, v_g'(t)\\
&= n\, \lambda^{n} \exp_g(\lambda;t)+\lambda^n\, \left(
 \lambda\, \varphi(t) \exp_g(\lambda;t) + \exp_g(\lambda;t)\right)\\
 &=(n+1)\, \lambda^{n} \exp_g(\lambda;t)+\lambda^{n+1}\, v(t),
\end{aligned}
\end{displaymath}
which finishes the proof.
\end{proof}

\section{The $g$-Wronskian and second order linear Stieltjes differential equations}

In this section we will define the concept of $g$-Wronskian and simplified $g$-Wronskian, and we will study their applications for second order linear Stieltjes differential equations. Observe that our definition of simplified $g$-Wronskian matches the definition of Wronskian in \cite[Définition~5.2.1]{Lariviere}. Nevertheless, our definition is more general as a consequence of having a broader definition of Stieljtes derivative that includes the points of $C_g$.

\begin{dfn}[$g$-Wronskian and simplified 
$g$-Wronskian] Given $y_1,\,y_2\in \mathcal{BC}_g^2([0,T];\mathbb{C})$, we define the 
\emph{$g$-Wronskian} as the map $W_g(y_1,y_2): [0,T]\to \mathbb C$ given by the expression
\begin{equation} \label{explicitgWrons}
W_g(y_1,y_2)(t)=
\begin{vmatrix}
y_1(t) && y_2(t) && (\Delta g(t))^2\\
(y_1)'_g(t) && (y_2)'_g(t) && -\Delta g(t)\\
(y_1)''_g(t) && (y_2)''_g(t) && 1\\
\end{vmatrix},\quad t\in [0,T].
\end{equation}
Explicitly, with the notation $\Delta g^2(t)=(\Delta g(t))^2$,
\begin{displaymath}
W_g(y_1,y_2)=
y_1 \, (y_2)'_g - y_2 \, (y_1)'_g  + [y_1\, (y_2)''_g- y_2\, (y_1)''_g]\, \Delta g  + [(y_1)'_g\, (y_2)''_g - (y_2)'_g\, (y_1)''_g]\, \Delta g^2.
\end{displaymath}

Similarly, given $y_1,y_2\in \mathcal{BC}_g^1([0,T];\mathbb{C})$, we define the \emph{simplified $g$-Wronskian}  as the function $\widetilde{W}_g(y_1,y_2): [0,T]\to \mathbb R$ given by the expression
\begin{displaymath}
\widetilde{W}_g(y_1,y_2)(t)=
\begin{vmatrix}
y_1(t) && y_2(t) \\
(y_1)'_g(t) && (y_2)'_g(t) 
\end{vmatrix},\quad t\in [0,T].
\end{displaymath}
Explicitly,
\begin{displaymath}
\widetilde{W}_g(y_1,y_2)(t)=y_1 \, (y_2)'_g - y_2 \, (y_1)'_g.
\end{displaymath}
\end{dfn}

\begin{rem} Observe that when condition \eqref{condDeltaderivable} is satisfied we can rewrite the $g$-Wronskian as
	\[W_g(y_1,y_2)^*=
	\begin{vmatrix}
		y_1 && y_2 && \left(\Delta g^*\right)^2\\
		(y_1)'_g && (y_2)'_g && (\Delta^2 g^*)_g'\\
		(y_1)''_g && (y_2)''_g && 1\\
	\end{vmatrix}.
	\]
	Under this condition, and noting that $(\Delta^2 g^*)''_g=
	(-\Delta g^*)'_g=\irchi_{D_g}^*$, we can further rewrite it as
	\[
	W_g(y_1,y_2)^*=(1-\irchi_{D_g}^*)
	\begin{vmatrix}
		y_1 && y_2 && 0\\
		(y_1)'_g && (y_2)'_g && 0\\
		(y_1)''_g && (y_2)''_g && 1\\
	\end{vmatrix}
	+\irchi_{D_g}^*
	\begin{vmatrix}
		y_1 && y_2&& \Delta^2 g^*\\
		(y_1)'_g && (y_2)'_g && (\Delta^2g^*)'_g\\
		(y_1)''_g && (y_2)''_g && (\Delta^2 g^*)''_g\\
	\end{vmatrix},\]
	or equivalently,
	\[
	W_g(y_1,y_2)^*=\begin{vmatrix}
		y_1 && y_2 && \Delta^2 g^*&&0\\
		(y_1)'_g && (y_2)'_g && (\Delta^2g^*)'_g&&0\\
		(y_1)''_g && (y_2)''_g && (\Delta^2 g^*)''_g&&1\\
		0&&0&&-(1-\irchi_{D_g}^*)&&\irchi_{D_g}^*
	\end{vmatrix}.
	\]
\end{rem}

\begin{rem}\label{WtildeACg}
	Observe that, given the order of derivation necessary in each of the definitions, we require different regularity on the functions for the definition of the $g$-Wronskian and the simiplified $g$-Wronskian. 
	Furthermore, we also need to note that $W_g(y_1,y_2)$ does not belong, in 
	general, to the space $\mathcal{BC}_g([0,T];\mathbb{C})$. However, thanks to the fact that $\mathcal{BC}_g^1([0,T];\mathbb{C})\subset 
	\mathcal{AC}_g([0,T],\mathbb{C})$ and \cite[Proposition 5.4]{FriLo17}, we have that
	$\widetilde{W}_g(y_1,y_2)\in 
	\mathcal{AC}_g([0,T];\mathbb{C})\subset \mathcal{BC}_g([0,T];\mathbb{C})$.
\end{rem}

It is easy to see that both the $g$-Wronskian and the simplified $g$-Wronskian yield the usual Wronskian of two functions when we consider $g=\operatorname{Id}$, that is, when the Stieltjes derivative coincides with the usual derivative. This can also be noted if we rewrite $W_g(y_1,y_2)$ as
\begin{displaymath}	
W_g(y_1,y_2)=\begin{vmatrix}
y_1+(y_1)'_g\Delta g &&& y_2+(y_2)'_g\Delta g\\
(y_1)'_g+(y_1)''_g\Delta g &&& (y_2)'_g+(y_2)''_g\Delta g
\end{vmatrix}.
\end{displaymath}
Bearing in mind Remark~\ref{remNgderivative}, we have that $f'_g(t)\Delta g(t)=f(t^+)-f(t)$ for $t\in D_g$, which means
\begin{displaymath}
W_g(y_1,y_2)(t)=\left\{
\begin{array}{ll}
y_1(t) (y_2)'_g(t)-y_2(t)(y_1)'_g(t),& \quad t\not\in [0,T]\cap D_g,\\
y_1(t^+) (y_2)'_g(t^+)-y_2(t^+)(y_1)'_g(t^+),&\quad  t\in [0,T]\cap D_g,
\end{array}\right.
\end{displaymath}
so, since $y_1,\,y_2\in \mathcal{BC}_g^2([0,T],\mathbb{C})$, Proposition~\ref{proreg} allows us to rewrite $W_g(y_1,y_2)$ simply as
\begin{equation}\label{gWronskianwithlimits}
W_g(y_1,y_2)(t)=y_1(t^+) (y_2)'_g(t^+)-y_2(t^+)(y_1)'_g(t^+),\quad
t\in [0,T],
\end{equation}
which shows that $W_g(y_1,y_2)(t)=\widetilde W_g(y_1,y_2)(t^+)$, $t\in[0,T]$. Furthermore, \eqref{gWronskianwithlimits} leads to the following result.
	
\begin{lem}\label{Wrightcont}
	Given $y_1,y_2\in \mathcal{BC}^2_g([0,T],\mathbb C)$,  
\[W_g(y_1,y_2)(t^+)=W_g(y_1,y_2)(t),\quad t\in[0,T).\]
As a consequence, $(W_g(y_1,y_2))'_g(t)=0$ for all $t\in[0,T)\cap D_g$.
	\end{lem}
\begin{proof}
	Taking into account equation \eqref{gWronskianwithlimits}, it is enough to check that 
	\[\lim_{s\to t^+}y_k(s^+)=y_k(t^+),\quad \lim_{s\to t^+}y_k'(s^+)=y_k'(t^+),\quad k=1,2.\]
	This follows directly from \cite[Corollary 4.1.9]{MonSlaTvr18} since the maps $y_k,y_k'$ are regulated as a direct consequence of Remark~\ref{remNgderivative} and Proposition~\ref{proreg}.
	\end{proof}

\begin{rem}
	Lemma~\ref{Wrightcont} ensures that $W_g(y_1,y_2)$ is continuous from the right at every $t\in[0,T)$. However, in general, we cannot ensure continuity at such points. Indeed, consider the maps
	\begin{displaymath}
		g(t)=\left\{
		\begin{array}{ll}
			t,\quad & t\le 0\\
			t+1,\quad & t>0,
		\end{array}
		\right.
		\quad
		y_1(t)=g|_{[-1,1]},\quad
		y_2(t)=\left\{
		\begin{array}{ll}
			e^{t+1},\quad & -1\le t\le 0,\\
			2e^{t+1},\quad & 0<t\leq1.
		\end{array}
		\right.
	\end{displaymath}
	It is possible to check that $y_1,y_2\in\mathcal{BC}^2_g([-1,1],\mathbb C)$ 
	and
	\begin{displaymath}
		(y_1)'_g(t)=1,\quad (y_1)''_g(t)=0,\quad (y_2)'_g(t)=(y_2)''_g(t)=y_2(t),\quad t\in[-1,1].
	\end{displaymath}
	Therefore, we can obtain the explicit expression of $W_g(y_1,y_2)$ from \eqref{explicitgWrons}:
	\begin{displaymath}
		W_g(y_1,y_2)(t)=\left\{
		\begin{array}{ll}
			e^{t+1}(t-1),\quad & -1\le t< 0,\\
			0,\quad & t=0,\\
			2te^{t+1},\quad & 0<t\leq1.
		\end{array}
		\right.
	\end{displaymath}
	In this case, we have that $W_g(y_1,y_2)(0^+)=W_g(y_1,y_2)(0)=0\not= -e=W_g(y_1,y_2)(0^-)$, which shows that $W_g(y_1,y_2)$ needs not be continuous at the points of $[-1,1]\cap D_g$. Furthermore, this also shows that $W_g(y_1,y_2)$ might not be $g$-continuous since it is not left-continuous at $0$, see Proposition~\ref{proreg}.
\end{rem}	

Similarly to the usual setting, our definition of Wronskian functions allows us to obtain a sufficient condition for two maps to be linearly independent in the corresponding space of functions, as presented in the next result.

\begin{lem}\label{lemwr}
	Let $y_1,y_2:[0,T]\to\mathbb C$. Then:
	\begin{itemize}
		\item If $y_1,y_2\in \mathcal{BC}_g^1([0,T];\mathbb{C})$ and there exists $t\in[0,T]$ such that $\widetilde{W}_g(y_1,y_2)(t)\ne 0$, then $y_1$ and $y_2$ are linearly independent. 
		\item If $y_1,y_2\in \mathcal{BC}_g^2([0,T];\mathbb{C})$ and there exists $t\in[0,T]$ such that $W_g(y_1,y_2)(t)\ne 0$, then $y_1$ and $y_2$ are linearly independent.
	\end{itemize}
\end{lem}
\begin{proof}
	First, suppose $y_1,y_2\in \mathcal{BC}_g^1([0,T];\mathbb{C})$ and there exists $t\in[0,T]$ such that $\widetilde{W}_g(y_1,y_2)(t)\ne 0$. Reasoning by contradiction, assume that $y_1$ and $y_2$ are linearly dependent. In that case, by the linearity of the $g$-derivative, the two first columns of the determinant that defines $\widetilde{W}_g(y_1,y_2)(t)$ must be linearly dependent for every $t\in[0,T]$, which implies that $\widetilde{W}_g(y_1,y_2)(t)=0$ for all $t\in[0,T]$, and this contradicts the hypotheses.
	
	Now, the proof for the case where $y_1,y_2\in \mathcal{BC}_g^2([0,T];\mathbb{C})$ and there exists $t\in[0,T]$ such that $W_g(y_1,y_2)(t)\ne 0$ is analogous and we omit it.
\end{proof}

The next result provides a partial converse to 	Lemma~\ref{lemwr} for functions in $\mathcal{BC}_g^1([0,T];\mathbb{C})$.

\begin{lem}\label{lempc}Let  $y_1,y_2\in \mathcal{BC}_g^1([0,T];\mathbb{C})$ be such that $y_1(0)\,y_2(0)\ne 0$, $y_1(t)\,y_2(t)\ne0$ for $g\mbox{-a.a.}\, t\in [0,T)$, and 
$(y_1)_g'/y_1\in \mathcal{L}_g^1([0,T),\mathbb C)$. If $\widetilde W_g(y_1,y_2)(t)= 0$, $g\mbox{-a.a.}\, t\in[0,T)$, then $y_1$ and $y_2$ are 
linearly dependent.
\end{lem}
\begin{proof} From the hypotheses, we know that there exists a $g$-measurable set, $N\subset[0,T)$,  such that $\mu_g(N)=0$ and $y_1(t)\,y_2(t)\ne 0$ for all $t\in[0,T)\setminus N$. Observe that, since $\widetilde{W}_g(y_1,y_2)(t)= 0$ $g\mbox{-a.a.}\, t\in[0,T)$, we must have that 
\begin{displaymath}
\frac{(y_1)'_g(t)}{y_1(t)}=\frac{(y_2)'_g(t)}{y_2(t)},\quad t\in [0,T)\setminus N.
\end{displaymath}
Define 
\begin{displaymath}
f:t\in [0,T)\rightarrow f(t)=\left\{
\begin{array}{ll}
 \displaystyle \frac{(y_1)'_g(t)}{y_1(t)}, & t\in [0,T)\setminus N,\\
0, & t\not\in N.
\end{array}
\right.
\end{displaymath}
Note that $f$ coincides with $\frac{(y_1)'_g}{y_1}$ on $[0,T)\setminus N$. Now, since $N$ is $g$-measurable and $\mu_g(N)=0$, it follows from  the hypotheses that $f\in\mathcal L^1_g([0,T),\mathbb C)$. Observe that $\frac{y_1}{y_1(0)}$ and $\frac{y_2}{y_2(0)}$ belong to $\mathcal{BC}_g^1([0,T];\mathbb{C})\subset\mathcal{AC}_g([0,T];\mathbb{C})$ and are solutions of the differential problem
\[u_g'(t)=f(t)\, u(t),\quad  g\mbox{-a.a. } t \in [0,T),\quad \quad u(0)=1,\]
which has a unique solution, ---cf. \cite[Theorem~4.6]{Ma21}. Therefore, $y_2(0)y_1-y_1(0)y_2=0$, which means that $y_1$ and $y_2$ are linearly dependent.
\end{proof}

\begin{rem}A result analogous to Lemma~\ref{lempc} can be stated for the case where $y_1,y_2\in\mathcal{BC}_g^2([0,T];\mathbb{C})$ and  $W_g(y_1,y_2)=0$, $\forall t\in [0,T]$.
	 Indeed, under these conditions, observing that $\widetilde W_g(y_1,y_2)$ is the principal minor determinant of $W_g(y_1,y_2)$ of order two, we have that, on $[0,T]$,
\begin{displaymath}
\widetilde{W}_g(y_1,y_2)= [y_2\, (y_1)''_g-y_1\, (y_2)''_g]\, \Delta g  + [(y_2)'_g\, (y_1)''_g-(y_1)'_g\, (y_2)''_g ]\, \Delta g^2,
\end{displaymath}
which, in particular, yields that $\widetilde{W}_g(y_1,y_2)(t)=0$, 
for $\forall t\in [0,T]\setminus D_g$. On the other hand, 
if $t\in [0,T]\cap D_g$, since $\widetilde{W}_g(y_1,y_2)$ is 
$g$-continuous, see Remark~\ref{WtildeACg}, Proposition~\ref{proreg} 
ensures that
\[\widetilde{W}_g(y_1,y_2)(t)=\lim_{s\to t^-} \widetilde{W}_g(y_1,y_2)(s)=\lim_{\substack{s\to t^-\\s\in [0,T]\setminus D_g}}\widetilde{W}_g(y_1,y_2)(s)=0,\]
so $\widetilde{W}_g(y_1,y_2)(t)=0$ $\forall t\in [0,T]$.
\end{rem}

So far, we have studied the properties of the $g$-Wronskian as a function. Now, we turn our attention to how this concept relates to Stieltjes differential equations. To that end, let us consider the following second order linear Stieltjes differential 
equation with $g$-continuous coefficients,
\begin{empheq}[left=\empheqlbrace]{align} 
& v_g''(t)+P(t)\, v_g'(t)+Q(t)\, v(t)=f(t),\quad  t \in [0,T],
\label{eq:secondordernh2A}\\
& v(0)=x_0,\;v_g'(0)=v_0, \label{eq:secondordernh2B}
\end{empheq}
where $x_0,\, v_0 \in \mathbb{C}$ and $f,\,P,\,Q\in \mathcal{BC}_g([0,T];\mathbb{C})$. We define the concept of solution in the following terms.

\begin{dfn} 
	A \emph{solution} of \eqref{eq:secondordernh2A} is a  function $v\in \mathcal{BC}_g^2([0,T];\mathbb{C})$ such that
	\[v_g''(t)+P(t)\, v_g'(t)+Q(t)\, v(t)=f(t),\quad  t \in [0,T].\]
	If, in addition, $v(0)=x_0$, $v'_g(0)=v_0$, then $v$ is a \emph{solution} of~\eqref{eq:secondordernh2A}-\eqref{eq:secondordernh2B}.
\end{dfn}

\begin{rem} \label{remexistence} Observe that~\eqref{eq:secondordernh2A}-\eqref{eq:secondordernh2B} is equivalent to the system 
	\begin{equation*}
		\left\{\begin{aligned}
			&x'_g(t)=y(t)&\quad  t \in [0,T],\\
			& y_g'(t)=-P(t) y(t)-Q(t)x(t)+f(t),&\quad  t \in [0,T],\\
			& x(0)=x_0,\ y(0)=v_0,&
		\end{aligned}\right.
	\end{equation*}
which we know to have a unique solution, 
---cf. \cite[Theorem 5.58]{MarquezTesis}. Thus, \eqref{eq:secondordernh2A}-\eqref{eq:secondordernh2B} has a unique solution.
\end{rem}

\begin{rem} 
	In general, we cannot ensure that a solution of \eqref{eq:secondordernh2A}-\eqref{eq:secondordernh2B}  belongs to space 
$\mathcal{BC}_g^{\infty}([0,T];\mathbb{C})$ even when we
consider $P,\,Q,\, f\in \mathcal{BC}_g^{\infty}([0,T];\mathbb{C})$ as the product of two functions in $\mathcal{BC}_g^\infty([0,T];\mathbb{C})$ 
is not necessarily a function in the space 
$\mathcal{BC}_g^\infty([0,T];\mathbb{C})$ 
(see \cite[Remark 3.16]{Fernandez2021} and \cite{Fernandez2022b}).
However, when  $P$ and $Q$ are constant, the regularity of the solution is determined by the regularity of $f$, see \cite{Fernandez2021}.
\end{rem}

In order to study the application of the variation of parameters method to obtain the solution of~\eqref{eq:secondordernh2A}-\eqref{eq:secondordernh2B}, we consider the homogeneous problem
\begin{empheq}[left=\empheqlbrace]{align} 
& v_g''(t)+P(t)\, v_g'(t)+Q(t)\, v(t)=0,\quad  t \in [0,T],
\label{eq:secondorderh2A}\\
& v(0)=x_0,\;v_g'(0)=v_0, \label{eq:secondorderh2B}
\end{empheq}
We have the following Lemma whose proof is straightforward from the linearity of the $g$-derivative.

\begin{lem}\label{homosol}
	 Let $y_1,y_2\in \mathcal{BC}_g^{2}([0,T];\mathbb{C})$ be two solutions of \eqref{eq:secondorderh2A} such that
\begin{equation}\label{Initialcondindep}
	y_1(0)(y_2)'_g(0)-y_2(0)(y_1)'_g(0)\neq 0.
\end{equation}
Then, the map $v=c_1 y_1+ c_2 y_2$ is the unique solution of~\eqref{eq:secondorderh2A}-\eqref{eq:secondorderh2B}, where
\[
c_1=\frac{(y_2)'_g(0)\, x_0-v_0\, y_2(0)}{y_1(0)\, (y_2)'_g(0)-y_2(0)\,(y_1)'_g(0)},\quad 
c_2=\frac{v_0\, y_1(0)-(y_1)'_g(0)\,x_0}{y_1(0)\, (y_2)'_g(0)-y_2(0)\,(y_1)'_g(0)}.
\]
\end{lem}

In the following lemma ---cf. \cite[Théorème 5.2.3]{Lariviere}--- we will make the relationship between 
$W_g(y_1,y_2)(t)$ and $\widetilde{W}_g(y_1,y_2)(t)$ explicit when 
$y_1(t)$ and $y_2(t)$ are solutions of the homogeneous equation.

\begin{lem} \label{GwronRel} 
	Let $y_1,y_2 \in \mathcal{BC}_g^{2}([0,T];\mathbb{C})$ be two solutions of \eqref{eq:secondorderh2A}. Then, 
\begin{equation}\label{WrelWtilde}
W_g(y_1,y_2)=\left(1-P\,\Delta g+Q\, \Delta g^2\right) 
\widetilde{W}_g(y_1,y_2).
\end{equation}
Moreover, if 
\begin{equation}\label{condPQ}
1-P(t)\,\Delta g(t)+Q(t)\, \Delta g(t)^2\neq 0,\quad t\in [0,T]\cap D_g,
\end{equation}
then, 
\begin{equation}\label{altexprWtilde}
	\widetilde{W}_g(y_1,y_2)=\widetilde{W}_g(y_1,y_2)(0)\,
	\exp_g(-P+Q\,\Delta g;\cdot).
\end{equation}
In particular, if $\widetilde{W}_g(y_1,y_2)(0)\neq 0$, then ${W}_g(y_1,y_2)(t)\ne 0$ for every $t\in[0,T]$ and the multiplicative inverse of  ${W}_g(y_1,y_2)$ is given by
\begin{equation}\label{multiplicativeinverseW}
W_g(y_1,y_2)^{-1}=\left[\widetilde{W}_g(y_1,y_2)(0)\, 
\left(1-P\,\Delta g+Q\, \Delta g^2\right) \right]^{-1}\, 
\exp_g\left(-\frac{-P+Q\,\Delta g}{1-P\,\Delta g+Q\, \Delta g^2} ;\cdot\right),
\end{equation}
which is a bounded function and continuous on $[0,T]\setminus D_g$.
\end{lem}

\begin{proof} Given that $y_1,y_2$ are solutions of \eqref{eq:secondorderh2A}, standard computations yield that
\begin{displaymath}
\begin{aligned}
y_1\, (y_2)''_g-y_2\, (y_1)''_g&=
-P\,\widetilde{W}_g(y_1,y_2), \\
(y_1)'_g\, (y_2)''_g-(y_2)'_g\, (y_1)''_g&=
Q \, \widetilde{W}_g(y_1,y_2),
\end{aligned}
\end{displaymath}
from which it is clear that \eqref{WrelWtilde} holds.

Assume that \eqref{condPQ} holds. In order to prove \eqref{altexprWtilde}, it is enough to show that $\widetilde{W}_g(y_1,y_2)$ is the unique solution of 
	\begin{equation*}\label{pw}
		\left\{\begin{aligned}
			u'_g(t)&=
			\left(-P(t)+Q(t)\, \Delta g(t)\right)u(t),\quad\quad 
			g\text{-a.a.}\, t\in [0,T), \\
			u(0)&=\widetilde{W}_g(y_1,y_2)(0).
		\end{aligned}\right.
\end{equation*}
This is because $\widetilde{W}_g(y_1,y_2)\in \mathcal{AC}_g([0,T];\mathbb{C})$ (see Remark~\ref{WtildeACg})  and, since \eqref{condPQ} holds, the unique solution in $\mathcal{AC}_g([0,T];\mathbb{C})$ of \eqref{pw} is $\widetilde{W}_g(y_1,y_2)(0)\exp_g\left(-P+Q\, \Delta g;t\right)$.

Clearly, $\widetilde{W}_g(y_1,y_2)$ satisfies the initial condition. Given the explicit expression of $\widetilde{W}_g(y_1,y_2)$, which is $g$-absolutely continuous on $[0,T]$, we can compute its derivative $g$-a.e. in $[0,T]$ by means of Proposition~\ref{PropStiDer} and Remark~\ref{f*equalsf}, which yield
\begin{displaymath}
	\begin{aligned}
		\left(\widetilde{W}_g(y_1,y_2)\right)'_g=&
		(y_1)'_g\, (y_2)'_g+y_1\, (y_2)''_g+
		(y_1)'_g\, (y_2)''_g\, \Delta g^*-(y_2)'_g\, (y_1)'_g-y_2\, (y_1)''_g-
		(y_2)'_g\, (y_1)''_g\, \Delta g^*
		\\
		=&\left(-P+Q\, \Delta g^*\right) \widetilde{W}_g(y_1,y_2),
	\end{aligned}
\end{displaymath}
so $	(\widetilde{W}_g(y_1,y_2))'_g=(-P+Q\, \Delta g) \widetilde{W}_g(y_1,y_2)$ $g$-a.e. in $[0,T]$. Now, thanks to~\eqref{medtneqts}, $\widetilde{W}_g(y_1,y_2)=\widetilde{W}_g(y_1,y_2)(0)\,
	\exp_g(-P+Q\,\Delta g;\cdot)$.

Finally, assume that \eqref{condPQ} holds and $\widetilde{W}_g(y_1,y_2)(0)\neq 0$. 
Given that \eqref{WrelWtilde} and\eqref{altexprWtilde} hold, it is clear that ${W}_g(y_1,y_2)(t)\ne 0$ for every $t\in[0,T]$. 
Observe that by Proposition~\ref{gexpprop}, we have that
\begin{displaymath}
	\exp_g(-P+Q\,\Delta g;\cdot)^{-1}=
	\exp_g\left(-\frac{-P+Q\,\Delta g}{1-P\,\Delta g+Q\, \Delta g^2} ;\cdot\right).
\end{displaymath}
Hence, given \eqref{WrelWtilde} and \eqref{altexprWtilde}, it follows that
\begin{align*}
	W_g(y_1,y_2)^{-1}&=\left(1-P\,\Delta g+Q\, \Delta g^2\right)^{-1} 
	\widetilde{W}_g(y_1,y_2)^{-1}\\
	&=\left[\widetilde{W}_g(y_1,y_2)(0)\, 
	\left(1-P\,\Delta g+Q\, \Delta g^2\right) \right]^{-1}\, 
	\exp_g\left(-\frac{-P+Q\,\Delta g}{1-P\,\Delta g+Q\, \Delta g^2} ;\cdot\right),
\end{align*}
which proves \eqref{multiplicativeinverseW}.  Hence, all that is left to do is show that $W_g(y_1,y_2)^{-1}$ is bounded and continuous at every $t\in[0,T]\setminus D_g$. 
Since $-P+Q\, \Delta g\in \mathcal{L}_g^1([0,T];\mathbb{C})$ and 
$1+\left(-P+Q\, \Delta g\right)\, \Delta g\neq 0$, we have, using the same 
arguments as in the proof of Corollary~\ref{lemtec1}, that 
$[1-P\, \Delta g+Q\, \Delta g^2]^{-1}$ is a bounded function, which ensures that $W_g(y_1,y_2)^{-1}$ is bounded. Finally, the continuity of $W_g(y_1,y_2)^{-1}$ on $[0,T]\setminus D_g$ is a direct consequence of Proposition~\ref{regulated}.	
\end{proof}

As a direct consequence of Lemma~\ref{GwronRel}, we have the following result ensuring the linear independence of two solutions of the homogeneous equation provided condition \eqref{Initialcondindep} is satisfied.

\begin{cor} \label{LinInd2} Assume that \eqref{condPQ} holds and let $y_1,y_2 \in \mathcal{BC}_g^{2}([0,T];\mathbb{C})$ be two solutions of \eqref{eq:secondorderh2A} such that \eqref{Initialcondindep} holds.
Then, $W_g(y_1,y_2)(t)\neq 0$ for all $t\in [0,T]$ and, in particular, 
 $y_1$ and $y_2$ are linearly independent.
\end{cor}

\begin{proof} Since $\widetilde{W}_g(y_1,y_2)(0)=y_1(0)(y_2)'_g(0)-y_2(0)(y_1)'_g(0)\neq 0$,
Lemma~\ref{GwronRel} guarantees that equation \eqref{altexprWtilde} holds, so, since \eqref{condPQ} holds and $\exp_g\left(-P+Q\, \Delta g;\cdot\right)\neq 0$, $W_g(y_1,y_2)(t)\neq 0$ for all $t\in [0,T]$. Now, the rest of the result follows from Lemma~\ref{lemwr}.
\end{proof}

\section{The variation of parameters method}

One of the most important applications of the $g$-Wronskian is its use in the variation of parameters method. This technique can be useful in two different ways: to find a linearly independent solution of the homogeneous equation~\eqref{eq:secondorderh2A} from a known solution of the same problem; or to obtain a particular solution of the nonhomogeneous problem~\eqref{eq:secondordernh2A} from two linearly independent solutions of the homogeneous problem.

\subsection{Finding a solution of  the homogeneous problem}

The idea behind the classical variation of parameters method is to find a linearly independent solution of~\eqref{eq:secondorderh2A} with respect to a known solution of~\eqref{eq:secondorderh2A}, say $y_1$, that can be  expressed as $y_2=\varphi\, y_1$ for a certain function $\varphi$. Here, we will show that a similar reasoning holds even in the context of Stieltjes calculus. 

	To that end, suppose that $y_1$ is a solution of equation \eqref{eq:secondorderh2A} and consider $y_2=\varphi\, y_1$ for 
	a bounded $g$-continuous function $\varphi$ twice $g$-differentiable. Let us deduce what kind of function $\varphi$ must be assuming that $y_2$ is a solution of \eqref{eq:secondorderh2A} which is linearly independent  from $y_1$.

First, under suitable conditions, we can use Propositions~\ref{PropStiDer} and \ref{diffDeltag*} to compute the first and second $g$-derivatives of $y_2$, which yield
\begin{align}
	(y_2)'_g&=\varphi_g'\, y_1+\varphi\, (y_1)'_g+\varphi_g'\, (y_1)'_g\, \Delta g^*,\label{y2firstderivative}\\
	(y_2)''_g&=\varphi_g''\, y_1+\varphi\, (y_1)''_g+\left(2-\chi_{D_g}^*\right) \, \varphi_g'\, (y_1)'_g+\varphi_g''\, (y_1)'_g\, \Delta g^*+\varphi_g'\, (y_1)''_g\, \Delta g^*.\nonumber
\end{align}
So, knowing that $y_1,y_2$ are solutions of equation \eqref{eq:secondorderh2A}, when substituting these expressions in \eqref{eq:secondorderh2A}, we obtain that $\varphi$ must satisfy the equation
\begin{displaymath}
	 \varphi_g''\,\left(y_1+(y_1)'_g\, \Delta g^*\right) +   \varphi_g'\, \left[P\, \left(y_1+(y_1)'_g\, \Delta g^* \right)
	+\left(2-\chi_{D_g}^*\right)\, (y_1)'_g+
	(y_1)''_g\, \Delta g^*\right]=0.
\end{displaymath}
Under suitable conditions, this is equivalent to requiring that $\varphi'_g=w$ for a certain function $w$ satisfying
\begin{align}
	w_g'&=-\left[P+\frac{\left(2-\chi_{D_g}^*\right)\, (y_1)'_g+(y_1)''_g\, \Delta g^*}{y_1+(y_1)'_g\, \Delta g^*}\right]\, w\nonumber\\
	&=-\bigg[P-Q\, \Delta g^*+\frac{(y_1)'_g}{y_1+(y_1)'_g\, \Delta g^*}
	\cdot \left[\left(2-\chi_{D_g}^*\right)
	-P\, \Delta g^*+Q\, (\Delta g^*)^2
	\right]\bigg]\, w,\label{eq:derw}
\end{align}
where we have used the fact that $y_1$ is a solution of \eqref{eq:secondorderh2A} to obtain the last equality. 

On the other hand, since $y_2$ is solution of \eqref{eq:secondorderh2A}, we know that $y_2\in \mathcal{BC}_g^2([0,T];\mathbb{C})$, so keeping in mind \eqref{y2firstderivative} and $w=\varphi'_g$, we can write $(y_2)'_g=\varphi (y_1)'_g+w(y_1+\Delta g^* (y_1)'_g)$,  which should belong to $\mathcal{BC}_g([0,T];\mathbb{C})$. In other words, we ought to have that 
\begin{displaymath}
	\eta:=w\left[y_1+\Delta g^*\, (y_1)'_g\right] \in \mathcal{BC}_g([0,T];\mathbb{C}).
\end{displaymath}
Observe that once again, under sufficient conditions, $\eta$ is $g$-differentiable on $[0,T]$ and,  by Propositions~\ref{PropStiDer} and \ref{diffDeltag*}, we have that
\begin{align*}
	\eta'_g=w&_g'\left[y_1+\Delta g^*\, (y_1)'_g\right]+
	w\left[(y_1)'_g-\chi_{D_g}^*(y_1)'_g+\Delta g^*\, (y_1)''_g-
	\Delta g^*\, (y_1)''_g\right]\\
	&+w_g'\left[(y_1)'_g-\chi_{D_g}^*(y_1)'_g+\Delta g^*\, (y_1)''_g-
	\Delta g^*\, (y_1)''_g\right]\Delta g^*
	\\=&w_g'\left[y_1+\Delta g^*\, (y_1)'_g\right]+
	w\left[(y_1)'_g-\chi_{D_g}^*(y_1)'_g\right]+
	w_g'\left[\Delta g^*\,(y_1)'_g-\Delta g^*\,(y_1)'_g\right]
	\\=&
	w_g'\left[y_1+\Delta g^*\, (y_1)'_g \right]+w\, \left(1-\chi_{D_g}^*\right)\, (y_1)'_g.
\end{align*}
Thus, keeping in mind equation~\eqref{eq:derw} and the definition of $\eta$,
\begin{align*}	
\eta'_g&=  -\left[P-Q\, \Delta g^*+\frac{(y_1)'_g}{y_1+(y_1)'_g\, \Delta g^*}
		\cdot \left[\left(2-\chi_{D_g}^*\right)
		-P\, \Delta g^*+Q\, (\Delta g^*)^2
		\right]\right]\, \eta \\ 
		&\ \ \ \  +\left(1-\chi_{D_g}^*\right)\, (y_1)'_g\, \frac{\eta}{y_1+(y_1)'_g\, \Delta g^*}\\
&=-\left[P-Q\, \Delta g^*
	+\frac{(y_1)'_g}{y_1+(y_1)'_g\,\Delta g^*} \left(1-P\,\Delta g^*+Q\, (\Delta g^*)^2\right)
	\right]\, \eta,
\end{align*}
which is a first order linear equation. This means that, under suitable conditions, we can consider
\begin{align*}
	\eta&=\exp_g\left(-P+Q\, \Delta g -\frac{(y_1)'_g}{y_1+(y_1)'_g\,\Delta g^*} \left(1-P\,\Delta g^*+Q\, (\Delta g^*)^2\right);\cdot \right)\\
	&=\exp_g\left(-P+Q\, \Delta g;\cdot\right) \, \exp_g\left(
	-\frac{(y_1)'_g}{y_1+(y_1)'_g\,\Delta g};\cdot\right),
\end{align*}
where the last equality is a consequence of the product formula in Proposition~\ref{gexpprop} and~\eqref{medtneqts}. Furthermore, the quotient rule in Proposition~\ref{gexpprop} ensures that, under the corresponding conditions,
\begin{displaymath}
	\left(\frac{1}{y_1}\right)'_g=-\left(\frac{(y_1)'_g}{y_1+(y_1)'_g\,\Delta g^*}\right)\,\frac{1}{y_1},
\end{displaymath}
so, given the uniqueness of solution of first order linear equations, it follows that 
\[\exp_g\left(-\frac{(y_1)'_g}{y_1+(y_1)'_g\,\Delta g};\cdot\right)=\frac{1}{y_1},\]
which means that $\eta$ can be expressed as
\begin{displaymath}
\eta=\frac{\exp_g\left(-P+Q\, \Delta g;\cdot\right)}{y_1},
\end{displaymath}
and, as a consequence,
\begin{displaymath}
w=\frac{\exp_g\left(-P+Q\, \Delta g;\cdot\right)}{y_1\left[y_1+ 
(y_1)'_g\, \Delta g^*\right]}.
\end{displaymath}
Now, the expression of $\varphi$ follows, yielding
\begin{equation} \label{eq:vvarpardef}
\varphi(t)=\int_{[0,t)} \frac{\exp_g\left(-P+Q\, \Delta g;s\right)}{y_1(s)\left[y_1+ (y_1)'_g(s)\, \Delta g(s)\right]}\, \operatorname{d}\mu_g(s).
\end{equation}

We must notice that in the formal deduction of 
the formula~\eqref{eq:vvarpardef} that we just derived, it has been necessary to obtain the $g$-derivative of 
$\Delta g^*$ (see Proposition \ref{diffDeltag*}). In the following result 
we will provide the explicit conditions under which we can obtain a 
linearly independent solution of the homogeneous 
problem~\eqref{eq:secondorderh2A} from a known solution of 
the same problem. We will see that it is not necessary to assume additional hypotheses that guarantee the $g$-derivability of $\Delta g^*$ in 
order to ensure that $y_2=\varphi\, y_1$ is a solution of \eqref{eq:secondorderh2A} which is linearly independent from $y_1$.

\begin{thm}\label{solhomoconog} 
	Assume that 
	\eqref{condPQ} holds.
	If $y_1$ is a solution of \eqref{eq:secondorderh2A} such that 
	$y_1(t)\neq 0$ for every $t \in [0,T]$; the map
	\begin{displaymath}
		t\in [0,T]\rightarrow \frac{1}{y_1(t)+(y_1)'_g(t)\, \Delta g(t)} 
	\end{displaymath}
	is well-defined and bounded; and 
	\begin{displaymath}
		t\in [0,T]\rightarrow\frac{1}{y_1(t)}
	\end{displaymath}
	belongs to $\mathcal{BC}_g([0,T];\mathbb{C})$, then the map $\varphi:[0,T]\to\mathbb C$ given by
	\begin{displaymath}
		\varphi(t)=\int_{[0,t)}\frac{ \exp_g\left(-P+Q\, \Delta g;s\right)}{y_1(t)\left[y_1(s)+(y_1)'_g(s)\, \Delta g(s)\right]} \, \operatorname{d}\mu_g(s)
		\in \mathcal{AC}_g([0,T];\mathbb{C}),
	\end{displaymath}
	is well-defined; and the map $y_2:=\varphi\, y_1$ is a solution of \eqref{eq:secondorderh2A}. Furthermore, $y_1$ and $y_2$ are linearly independent.
\end{thm}
\begin{proof} On the one hand, thanks to \cite[Lemma~2.14]{Fernandez2022b}, we have that the function
$t\in [0,T]\rightarrow  \exp_g\left(-P+Q\, \Delta g;t\right)/y_1$ belongs to $\mathcal{BC}_g([0,T];\mathbb{C})$, so thanks to Corollary~\ref{lemtec2}, $y_2=\varphi\, y_1\in \mathcal{BC}_g^1([0,T];\mathbb{C})$ and
\begin{displaymath}
\begin{aligned}
(y_2)'_g=& (y_1)'_g\, \varphi+\frac{1}{y_1}\, \exp_g\left(-P+Q\, \Delta g;\cdot\right). 
\end{aligned}
\end{displaymath}
Now, thanks to Proposition~\ref{PropStiDer},
\begin{displaymath}
\begin{aligned}
(y_2)''_g=& (y_1)''_g\, \varphi+   \frac{(y_1)'_g\, \exp_g\left(-P+Q\, \Delta g;\cdot\right)}{y_1\left[y_1+(y_1)_g'\,\Delta g^*\right]} + \frac{(y_1)''_g\, \exp_g\left(-P+Q\, \Delta g;\cdot\right)\, \Delta g^*}{y_1\left[y_1+(y_1)_g'\,\Delta g^*\right]} 
\\
&+\frac{y_1\,\left( -P+Q\, \Delta g^* \right)\, \exp_g\left(-P+Q\, \Delta g;\cdot\right)}{
y_1\left[y_1+(y_1)_g'\,\Delta g^*\right]}-\frac{(y_1)'_g\, \exp_g\left(-P+Q\, \Delta g;\cdot\right)}{y_1\left[y_1+(y_1)_g'\,\Delta g^*\right]}.
\end{aligned}
\end{displaymath}
Hence,
\begin{displaymath}
(y_2)''_g= (y_1)''_g\, \varphi+ 
\frac{\exp_g\left(-P+Q\, \Delta g;\cdot\right)}{y_1\left[y_1+(y_1)_g'\,\Delta g^*\right]}
\, \left[(y_1)''_g\, \Delta g^*+y_1\,\left( -P+Q\, \Delta g^* \right)\right].
\end{displaymath}
If we use the fact that $(y_1)''_g=-P \, (y_1)'_g-Q\, y_1$, we derive that
\begin{displaymath}
\begin{aligned}
 (y_1)''_g\, \Delta g^*+y_1\,\left( -P+Q\, \Delta g^* \right)= & -P \, (y_1)'_g\, \Delta g^*-Q\, y_1\, \Delta g^*+y_1\,\left( -P+Q\, \Delta g^* \right) \\
=&-P\left[y_1+(y_1)_g'\,\Delta g^*\right].
\end{aligned}
\end{displaymath}
As a consequence,
\begin{displaymath}
(y_2)''_g=(y_1)''_g\, \varphi-\frac{P}{y_1} \, \exp_g\left(-P+Q\, \Delta g;\cdot\right).
\end{displaymath}
Observe that $(y_2)''_g\in \mathcal{BC}_g([0,T];\mathbb{C})$, so $y_2\in \mathcal{BC}_g^2([0,T];\mathbb{C})$.
 Let us now see that $y_2$ satisfies equation~\eqref{eq:secondorderh2A}. Indeed,
\begin{displaymath}
\begin{aligned}
& (y_2)''_g+P\, (y_2)'_g+Q\, y_2 \\ = & (y_1)''_g\, \varphi-\frac{P}{y_1} \, \exp_g\left(-P+Q\, \Delta g;\cdot\right) +P\,\left[(y_1)'_g\, \varphi+\frac{1}{y_1}\, \exp_g\left(-P+Q\, \Delta g;\cdot\right) \right]+Q\, \varphi\, y_1 \\
=&\varphi\, \left[ (y_1)''_g+P\, (y_1)'_g+Q\, y_1\right]=0,
\end{aligned}
\end{displaymath}
since $y_1$ is a solution of~\eqref{eq:secondorderh2A}.

Finally,  we  study the linear independence by means of Lemma~\ref{lemwr}. Observe that
\begin{displaymath}
\begin{aligned}
\widetilde{W}_g(y_1,y_2)=& y_1\, (y_2)'_g-y_2\, (y_1)'_g\\
=& y_1\, \left[ (y_1)'_g\, \varphi+y_1\, \varphi_g'+(y_1)'_g\,\varphi_g'\, \Delta g^*\right]-y_1\,\varphi\, (y_1)'_g\\
=& y_1\, \varphi_g'\left[y_1+(y_1)'_g\, \Delta g^*\right]\\
=&\exp_g\left(-P+Q\, \Delta g;t\right)\neq 0,
\end{aligned}
\end{displaymath}
which is enough to guarantee  that $y_1$ and $y_2$ are linearly independent on $\mathcal{BC}_g([0,T];\mathbb{C})$.
\end{proof}

\begin{rem}
	Similarly to Corollary~\ref{lemtec2}, Theorem~\ref{solhomoconog} presents some technical conditions on the map $y_1$, which can be ignored if we can guarantee that $y_1$ is bounded away from zero. Indeed, if that is the case, it is clear that $y_1(t)\not=0$, $t\in[0,T]$, while Remark~\ref{remboundedfromzero} ensures that the map $(y_1+(y_1)'_g\Delta g)^{-1}$ is well-defined and bounded. Following a similar reasoning, we can also see that $y_1^{-1}$ is well-defined and bounded, while the fact that $y_1$ is $g$-continuous on $[0,T]$ and $y_1(t)\not=0$, $t\in[0,T]$, are enough to ensure that $y_1^{-1}$ is $g$-continuous on $[0,T]$.
\end{rem}

\subsection{Obtaining a particular solution of the nonhomogeneous problem}

In a similar way to the classical case, to obtain the solution of the nonhomogeneous equation we need to find a particular solution. In order to achieve this, we can use the method of variation of parameters, that is, we will look for a particular solution of the form
\begin{equation}\label{eq:particular1}
v_p=c_1\, y_1+c_2\, y_2.
\end{equation}
where $c_1,\, c_2\in \mathcal{BC}_g([0,T];\mathbb{C})$ are $g$-differentiable at all points of $[0,T]$. In that case, thanks to 
Proposition~\ref{PropStiDer}, we have
\begin{equation*}
(v_p)'_g=
\displaystyle
c_1\, (y_1)'_g+c_2\, (y_2)'_g  + (c_1)'_g \, y_1+ (c_2)'_g\, y_2
 + (c_1)'_g\, (y_1)'_g\,\Delta g^*+(c_2)'_g \, (y_2)'_g\, \Delta g^*.
\end{equation*}
In order to avoid that second derivatives appear on the unknowns $c_1$ and $c_2$ when computing the second derivative of $v_p$ we will assume that
\begin{equation} \label{eq:nohomo1}
(c_1)'_g \, y_1+ (c_2)'_g\, y_2 
+ (c_1)'_g\, (y_1)'_g\,\Delta g^*+(c_2)'_g \, (y_2)'_g\,\Delta g^*=0.
\end{equation}
Thus,
\begin{equation*}
(v_p)''_g=\displaystyle
(c_1)'_g\, (y_1)'_g+c_1\, (y_1)''_g+ (c_2)'_g\, (y_2)'_g+ c_2\, (y_2)''_g +(c_1)'_g\, (y_1)''_g\, \Delta g^*+(c_2)'_g\, (y_2)''_g\, \Delta g^*.
\end{equation*}
Substituting $v_p$ in (\ref{eq:secondordernh2A}) we obtain that
\begin{equation} \label{eq:nohomo2}
(c_1)'_g\, (y_1)'_g+(c_2)'_g\, (y_2)'_g+ (c_1)'_g\, (y_1)''_g\, \Delta g^*+(c_2)'_g\, (y_2)''_g\, \Delta g^*=
f.
\end{equation}
Taking into account equations (\ref{eq:nohomo1}) and (\ref{eq:nohomo2}), we need 
$(c_1)'_g$ and $(c_2)'_g$ to satisfy the linear system
\begin{equation*}
\left( \begin{array}{rr}
y_1+(y_1)'_g\, \Delta g^* & y_2+(y_2)'_g\, \Delta g^* \\
(y_1)'_g+(y_1)''_g\, \Delta g^* & (y_2)'_g+(y_2)''_g\, \Delta g^*
\end{array}\right) \left(\begin{array}{c}
(c_1)'_g\\
(c_2)'_g
\end{array}\right)=\left(\begin{array}{c}
0\\
f
\end{array}\right).
\end{equation*}
The determinant of the matrix associated to the previous system coincides with the
$g$-Wronskian evaluated at $t^*$, so, if we can ensure that $W_g(y_1,y_2)(t)\neq 0$, $t\in [0,T]$, then,
\begin{equation}\label{eq:gwronk1}
(c_1)'_g= \frac{-(y_2+(y_2)'_g\, \Delta g^*)\, f}{W_g(y_1,y_2)^*},\quad 
(c_2)'_g= \frac{(y_1+(y_1)'_g\, \Delta g^*)\, f}{W_g(y_1,y_2)^*}.
\end{equation}
In the next result we will see that the functions
\begin{align}
c_1(t)&= \int_{[0,t)} \frac{-(y_2(s)+(y_2)'_g(s)\, \Delta g(s))\, f(s)}{W_g(y_1,y_2)(s)}\, \operatorname{d}\mu_g(s), \vspace{0.2cm}\label{eq:defc1}\\
c_2(t)&= \int_{[0,t)} \frac{(y_1(s)+(y_1)'_g(s)\, \Delta g(s))\, f(s)}{W_g(y_1,y_2)(s)}\, \operatorname{d}\mu_g(s),\label{eq:defc2}
\end{align}
are such that~\eqref{eq:particular1} is a particular solution of equation~\eqref{eq:secondordernh2A}. We can now state and proof 
the following Theorem for the solutions to the initial value problem~\eqref{eq:secondordernh2A}-\eqref{eq:secondordernh2B} using 
the definition of $g$-Wronskian. This generalizes \cite[Théorème~5.2.4]{Lariviere} where they considered constant coefficients.

\begin{thm} \label{solnohomoconog} Assume that~\eqref{condPQ} holds. Let $y_1,\, y_2\in \mathcal{BC}_g^2([0,T];\mathbb{C})$ two solutions of~\eqref{eq:secondorderh2A} such that $W_g(y_1,y_2)(t)\neq 0$ for 
all $t\in [0,T]$. Then,
\begin{equation}\label{spnh}
\begin{aligned}
v_p(t)&= y_1(t)\, \int_{[0,t)} \frac{-(y_2(s)+(y_2)'_g(s)\, \Delta g(s))\, f(s)}{W_g(y_1,y_2)(s)}\, \operatorname{d}\mu_g(s)\\
&\ \ \ \ +y_2(t)\, \int_{[0,t)} \frac{(y_1(s)+(y_1)'_g(s)\, \Delta g(s))\, f(s)}{W_g(y_1,y_2)(s)}\, \operatorname{d}\mu_g(s)
\end{aligned}
\end{equation}
is a particular solution to the nonhomogeneous 
equation~\eqref{eq:secondordernh2A}. Moreover, $v(t)=v_p(t)+v_h(t)$ is the solution of the initial value 
problem~\eqref{eq:secondordernh2A}-\eqref{eq:secondordernh2B}, where
\begin{displaymath} 
v_h(t)= \left(\frac{x_0\, (y_2)'_g(0)-v_0\, y_2(0)}{W_g(y_1,y_2)(0)}\right) \, y_1(t)+
\left(\frac{v_0\, y_1(0)-x_0\, (y_1)'_g(0)}{W_g(y_1,y_2)(0)}\right) \, y_2(t)
\end{displaymath}
is the solution of the initial value problem~\eqref{eq:secondorderh2A}-\eqref{eq:secondorderh2B}. 
\end{thm}

\begin{proof} On the one hand, thanks to Lemma~\ref{GwronRel},
\begin{displaymath}
W_g(y_1,y_2)^{-1}=
\left[\widetilde{W}_g(y_1,y_2)(0) \left(1-P\,\Delta g+Q\, \Delta g^2\right) \right]^{-1} \cdot \exp_g\left(-\frac{-P+Q\,\Delta g}{1-P\,\Delta g+Q\, \Delta g^2} ;\cdot\right)
\end{displaymath}
is a bounded function and continuous on $[0,T]\setminus D_g$. Hence, the functions
\[
\frac{-(y_2+(y_2)'_g\, \Delta g)\, f}{W_g(y_1,y_2)},\quad\quad
\frac{(y_1+(y_1)'_g\, \Delta g)\, f}{W_g(y_1,y_2)},
\]
satisfy the hypotheses of Lemma~\ref{dercont} and, therefore, $c_1$ and $c_2$, given by~\eqref{eq:defc1} and \eqref{eq:defc2} respectively, are $g$-absolutely continuous on $[0,T]$ and satisfy the equations in \eqref{eq:gwronk1}  for all $t\in [0,T]$. 

 Let us now check that $v_p\in \mathcal{BC}_g^2([0,T];\mathbb{C})$. Indeed, first observe that $v_p\in \mathcal{BC}_g([0,T]; \mathbb{C})$. Now, thanks to Proposition~\ref{PropStiDer}, we have that
\begin{displaymath}
\begin{aligned}
(v_p)'_g=& (y_1)'_g c_1+ (y_2)'_g c_2 - y_1 \frac{(y_2+(y_2)'_g\, \Delta g^*)\, f}{W_g(y_1,y_2)^*}+y_2 \frac{(y_1+(y_1)'_g\, \Delta g^*)\, f}{W_g(y_1,y_2)^*}\\&
- (y_1)'_g \frac{(y_2+(y_2)'_g\, \Delta g^*)\, f}{W_g(y_1,y_2)^*}\Delta g^*+(y_2)'_g \frac{(y_1+(y_1)'_g\, \Delta g^*)\, f}{W_g(y_1,y_2)^*}\Delta g^*\\= & (y_1)'_g c_1+ (y_2)'_g c_2.
\end{aligned}
\end{displaymath}
Therefore,
\begin{equation}\label{dspnh}
\begin{aligned}
(v_p)'_g(t)=& (y_1)'_g(t) \int_{[0,t)} \frac{-(y_2+(y_2)'_g\, \Delta g)\, f}{W_g(y_1,y_2)}\, \operatorname{d}\mu_g \\
& + (y_2)'_g(t) \int_{[0,t)} \frac{(y_1+(y_1)'_g\, \Delta g)\, f}{W_g(y_1,y_2)}\, \operatorname{d}\mu_g, 
\end{aligned}
\end{equation}
for all $t\in [0,T]$, so $(v_p)'_g\in 
\mathcal{BC}_g([0,T];\mathbb{C})$. Reasoning in a similar way,
\begin{displaymath}
\begin{aligned}
(v_p)''_g(t)&= (y_1)''_g(t) \int_{[0,t)} \frac{-(y_2+(y_2)'_g\, \Delta g)\, f}{W_g(y_1,y_2)}\, \operatorname{d}\mu_g \\
&\ \ \ \ + (y_2)''_g(t) \int_{[0,t)} \frac{(y_1+(y_1)'_g\, \Delta g)\, f}{W_g(y_1,y_2)}\, \operatorname{d}\mu_g+
f(t)\\
&=(y_1)''_g(t) c_1(t)+ (y_2)''_g(t)c_2(t)+f(t),
\end{aligned}
\end{displaymath}
for all $t\in [0,T]$. Whence, $(v_p)''_g\in \mathcal{BC}_g([0,T];\mathbb{C})$, which  shows that $v_p\in \mathcal{BC}^2_g([0,T];\mathbb{C})$.

Observe that it is immediate from equations \eqref{spnh} and \eqref{dspnh} that $v_p(0)=(v_p)'(0)=0$. 
Furthermore, $v_p$ satisfies, indeed, \eqref{eq:secondordernh2A} since
\begin{align*}
	(v_p)_g''+P (v_p)_g'+Q v_p&=(y_1)''_g c_1+ (y_2)''_gc_2+f++P )(y_1)'_g c_1+ (y_2)'_gc_2)+Q (y_1 c_1+ y_2c_2)\\
	&=c_1((y_1)''_g+P(y_1)'_g+Qy_1)+c_2((y_2)''_g+P(y_2)'_g+Qy_2)+f=f.
\end{align*}
This proves that $v_p$ is a solution of \eqref{eq:secondordernh2A}. Now, the rest of the result follows from Lemma~\ref{homosol}.
\end{proof}

\section{Applications} 

In this final section, we illustrate the theory developed above with three examples: a second order linear Stieltjes differential equation with constant coefficients, an example with variable coefficients and the one-dimensional linear Helmholtz equation with piecewise-constant coefficients.
	
\subsection{Second order linear Stieltjes differential equations with constant coefficients}

In \cite{Fernandez2021}, the authors obtained the solution of the second-order linear problem
\begin{displaymath}
	\left\{\begin{aligned}
		& v_g''(t)+ P\, v_g'(t)+Q\,v(t)=f(t),\quad t \in [0,T], \\
		& v(0)=x_0,\;v_g'(0)=v_0,
	\end{aligned}\right.\end{displaymath}
where $P,Q,x_0,v_0\in{\mathbb{C}}$, $f\in 
\mathcal{BC}_g([0,T];{\mathbb{C}})$. The approach that they took  was a factorization approach. Namely, considering $\lambda_1,\lambda_2\in\mathbb C$ such that for $k=1,2$,
\begin{align}
	&\lambda_k^2+P\lambda_k+Q=0\label{condlambda1}\\
	&1+\lambda_k \, \Delta g(t)\neq 0,\quad  t\in[0,T]\cap D_g,\label{condlambda2}
\end{align}
 they showed that a solution of the second order problem, $v$, must satisfy 
\begin{equation*}\left\{\begin{array}{l}
		v'_g(t)= \lambda_2\,v(t)+v_1(t)\quad t \in [0,T],\\
		v(0)= x_0,
	\end{array}\right.\end{equation*}
where $v_1=v_g'-\lambda_2\,v$ must satisfy
\begin{equation*}\left\{\begin{array}{l}
		(v_1)'_g(t)= \lambda_1v_1(t)+f(t),\quad  t \in [0,T],\\
		v(0)= v_0-\lambda_2x_0.
	\end{array}\right.\end{equation*}
By doing so, in \cite[Proposition 4.12]{Fernandez2021} they obtained the explicit expression of the solution, which is given by
\begin{equation} \label{general}
	\begin{aligned}
		v(t)&=x_0\exp_g(\lambda_2;t) +
		\left( v_0-\lambda_2x_0\right) \,\exp_g(\lambda_2;t)\, \int_{[0,t)} \frac{\exp_g
			\left( \frac{\lambda_1-\lambda_2}{1+\lambda_2\, \Delta g};\cdot\right)}{1+\lambda_2 \Delta g}\,  \, \operatorname{d} \mu_g\\
		&\ \  +\exp_g(\lambda_2;t)\, \int_{[0,t)} \frac{\exp_g
			\left( \frac{\lambda_1-\lambda_2}{1+\lambda_2\, \Delta g};\cdot\right)}{1+\lambda_2 \Delta g}\,  \cdot
		\bigg(\int_{[0,s)} \exp_g(\lambda_1;r)^{-1}\, \frac{f(r)}{1+\lambda_1 \Delta g(r)}\, \operatorname{d} \mu_g(r)\bigg)\, \operatorname{d} \mu_g.
\end{aligned}\end{equation}

In this section we  present an alternative method of obtaining the general solution given by~\eqref{general} 
based on the variation of parameters method. This technique will be especially useful when equation $x^2+Px+Q=0$ only has a single real root $\lambda=-P/2$. We have the following result.

\begin{cor} \label{thmhomo} Let $\lambda_1,\lambda_2\in \mathbb{C}$ be such that \eqref{condlambda1} and \eqref{condlambda2} hold. Then:
\begin{itemize}
\item If $\lambda_1 \neq \lambda_2$, then $y_k=\exp_g(\lambda_k;\cdot)$, $k=1,2$, are two linearly independent solutions of~\eqref{eq:secondorderh2A} which belong to $\mathcal{BC}_g^{\infty}([0,T];\mathbb{C})$.
\item If $\lambda_1=\lambda_2:=\lambda$, then $y_1=\exp_g(\lambda;\cdot)$ and 
\begin{equation} \label{sol2homo}
y_2(t)=\left[\int_{[0,t)} \frac{1}{1+\lambda \Delta g} \, \operatorname{d}\mu_g\right] \exp_g(\lambda;t),\quad t\in[0,T],
\end{equation}
are two linearly independent solutions of~\eqref{eq:secondorderh2A} which belong to $\mathcal{BC}_g^{\infty}([0,T];
\mathbb{C})$.
\end{itemize}
\end{cor}

\begin{proof} The first case is straightforward from Corollary~\ref{LinInd2} 
and Remark~\ref{expCinfty}. Just observe that 
\begin{displaymath}
y_1(0)(y_2)'_g(0)-y_2(0)(y_1)'_g(0)=\lambda_2-\lambda_1 \neq 0.
\end{displaymath}

Let us now consider the second case, namely, $\lambda_1=\lambda_2=-P/2:=\lambda$.  In that case, we know that $y_1(t)=\exp_g(\lambda;t)$ is a  solution of \eqref{eq:secondorderh2A}  and it belongs to $\mathcal{BC}_g^{\infty}([0,T];
\mathbb{C})$, see Remark~\ref{expCinfty}. We shall show that $y_2$ is the linearly independent solution given by Theorem~\ref{solhomoconog}.

First, observe that condition \eqref{condPQ} is satisfied since 
\begin{displaymath}
	1-P\,\Delta g(t)+Q\, \Delta g(t)^2=(1+\lambda\, \Delta g(t))^2\neq 0,\quad t \in [0,T].
\end{displaymath}
Furthermore, since $1+\lambda\, \Delta g(t)\neq 0$, 
$t\in [0,T]\cap D_g$, we have that $y_1(t)=\exp_g(\lambda;\cdot)\neq 0$, 
$t\in [0,T]$, and moreover, by Proposition~\ref{gexpprop}, $y_1^{-1}= \exp_g\left(-\frac{\lambda}{1+\lambda\, \Delta g};\cdot\right)$, so $y_1^{-1}\in\mathcal{BC}_g([0,T];\mathbb{C})$. 

On the other hand,
\begin{displaymath}
\frac{1}{y_1+(y_1)'_g\, \Delta g}=\frac{1}{y_1} \frac{1}{1+\lambda \, 
\Delta g}=\frac{\exp_g\left(-\frac{\lambda}{1+\lambda\, \Delta g};\cdot\right)}{1+\lambda \, 
\Delta g}
\end{displaymath}
is well-defined and bounded. Hence, the hypotheses of Theorem~\ref{solhomoconog} are satisfied which implies that $y_2=\varphi\, y_1$, with
\begin{displaymath}
\varphi(t)=\int_{[0,t)}\frac{ \exp_g\left(-P+Q\, \Delta g;\cdot\right)}{y_1\left[y_1+(y_1)'_g\, \Delta g\right]} \, \operatorname{d}\mu_g,
\end{displaymath}
is a solution of~\eqref{eq:secondorderh2A} and $y_1$ and $y_2$ are linearly independent. Let us show that $y_2$ can be expressed as \eqref{sol2homo}.

First, since $y_1=\exp_g(\lambda;\cdot)$, we have that $y_1'=\lambda y_1$ and, thus $y_1+y_1'\Delta g=y_1(1+\lambda\Delta g)$. On the other hand, since $P^2-4Q=0$ and $\lambda=-P/2$, we have that $-P+Q\Delta g=2\lambda+\lambda^2\Delta g$, so, using Proposition~\ref{gexpprop}, we see that
\begin{displaymath}
\begin{aligned}
\varphi(t)&=\int_{[0,t)}\frac{ \exp_g\left(2\, \lambda+\lambda^2\, \Delta g;\cdot\right)}{y_1\left[y_1+(y_1)'_g\, \Delta g\right]} \, \operatorname{d}\mu_g=\int_{[0,t)}\frac{ \exp_g(\lambda;\cdot)^{2}}{y_1^2\left[1+\lambda\, \Delta g\right]} \, \operatorname{d}\mu_g
=\int_{[0,t)} \frac{1}{1+\lambda \, \Delta g}\, \operatorname{d}\mu_g,
\end{aligned}
\end{displaymath}
so $y_2$ coincides with the expression in \eqref{sol2homo}.
Finally, Corollary~\ref{prodcinf} ensures that $y_2$ belongs to $\mathcal{BC}_g^{\infty}([0,T];\mathbb{C})$, which finishes the proof.
\end{proof}

Let us  analyze next the existence of a solution of the initial value problem~\eqref{eq:secondorderh2A}-\eqref{eq:secondorderh2B}. We have 
the following corollary.

\begin{cor}\label{constantcoeffsol}
	Let $\lambda_1,\lambda_2\in \mathbb{C}$ be such that \eqref{condlambda1} and \eqref{condlambda2} hold. Let $y_1,\,y_2\in \mathcal{BC}_g^2([0,T];\mathbb{C})$ be the two linearly independent solutions of~\eqref{eq:secondorderh2A} provided by \textup{Corollary~\ref{thmhomo}}. Then:
\begin{enumerate}
\item If $\lambda_1 \neq \lambda_2$, then, the map
\begin{displaymath}
\begin{aligned}
v(t)=&\frac{\lambda_2 x_0-v_0}{\lambda_2-\lambda_1}\exp_g(\lambda_1;t)+
\frac{v_0-\lambda_1 x_0}{\lambda_2-\lambda_1}\exp_g(\lambda_2;t)\\
&+\frac{1}{\lambda_1-\lambda_2} \exp_g(\lambda_1;t)\int_{[0,t)} \exp_g(\lambda_1;\cdot)^{-1}\, \frac{f}{1+\lambda_1 \Delta g}\, \operatorname{d} \mu_g\\
&-\frac{1}{\lambda_1-\lambda_2} \exp_g(\lambda_2;t)\int_{[0,t)} \exp_g(\lambda_2;\cdot)^{-1} \frac{f}{1+\lambda_2 \Delta g}
\, \operatorname{d} \mu_g.
\end{aligned}
\end{displaymath}
is the solution of the initial value problem~\eqref{eq:secondorderh2A}-\eqref{eq:secondorderh2B}.
\item If $\lambda_1=\lambda_2:=\lambda$, then, the map
\begin{displaymath}
\begin{aligned}
v(t)=&x_0 \, \exp_g(\lambda;t)+\left(v_0-\lambda x_0\right) \varphi(t) \exp_g(\lambda;t)\\
&- \exp_g(\lambda;t) \int_{[0,t)} 
\exp_g(\lambda;\cdot)^{-1}\,
\left[
\varphi\, \frac{f}{1+\lambda \Delta g}+
\frac{f\, \Delta g}{\left(1+\lambda \Delta g\right)^2}
 \right]\, \operatorname{d} \mu_g\\
 &+\varphi(t)\, \exp_g(\lambda;t) \, \int_{[0,t)} \exp_g(\lambda;\cdot)^{-1} \frac{f}{1+\lambda \Delta g}
 \, \operatorname{d} \mu_g,
\end{aligned}
\end{displaymath}
is the solution of the initial value problem~\eqref{eq:secondorderh2A}-\eqref{eq:secondorderh2B}. Where,
\begin{displaymath}
\varphi(t)=\int_{[0,t)} \frac{1}{1+\lambda \, \Delta g}\, \operatorname{d}\mu_g.
\end{displaymath}
\end{enumerate}
\end{cor}

\begin{proof} The result is a direct consequence of Theorem~\ref{solnohomoconog}. Indeed, let us consider the two cases separately.
\begin{enumerate}
\item If $\lambda_1\neq \lambda_2$,
we have that
\begin{displaymath}
\begin{aligned}
\frac{-(y_2(s)+(y_2)'_g(s)\, \Delta g(s))\, f(s)}{W_g(y_1,y_2)(s)}
=&+\frac{1}{\lambda_1-\lambda_2} \exp_g(\lambda_1;s)^{-1} \frac{f(s) }{1+\lambda_1\Delta g(s)},\\
\frac{(y_1(s)+(y_1)'_g(s)\, \Delta g(s))\, f(s)}{W_g(y_1,y_2)(s)}=&
-\frac{1}{\lambda_1-\lambda_2} \exp_g(\lambda_2;s)^{-1} \frac{f(s) }{1+\lambda_2\Delta g(s)}.
\end{aligned}
\end{displaymath}
Thus, Theorem~\ref{solnohomoconog} ensures that the solution is the map
\begin{equation}\label{solgB}
\begin{aligned}
v(t)=&\frac{\lambda_2 x_0-v_0}{\lambda_2-\lambda_1}\exp_g(\lambda_1;t)+
\frac{v_0-\lambda_1 x_0}{\lambda_2-\lambda_1}\exp_g(\lambda_2;t)\\
&+\frac{1}{\lambda_1-\lambda_2} \exp_g(\lambda_1;t)\int_{[0,t)} \exp_g(\lambda_1;s)^{-1}\, \frac{f(s)}{1+\lambda_1 \Delta g(s)}\, \operatorname{d} \mu_g(s)\\
&-\frac{1}{\lambda_1-\lambda_2} \exp_g(\lambda_2;t)\int_{[0,t)} \exp_g(\lambda_2;s)^{-1} \frac{f(s)}{1+\lambda_2 \Delta g(s)}
\, \operatorname{d} \mu_g(s).
\end{aligned}
\end{equation}
\item On the other hand, if $\lambda_1=\lambda_2$, we have
\begin{displaymath}
\begin{aligned}
\frac{-(y_2(s)+(y_2)'_g(s)\, \Delta g(s))\, f(s)}{W_g(y_1,y_2)(s)}=&
-\varphi(s)\,\exp_g(\lambda;s)^{-1}\, \frac{f(s)}{1+\lambda \Delta g(s)}\\
&-
\exp_g(\lambda;s)^{-1} \, \frac{f(s)\, \Delta g(s)}{\left(1+\lambda \Delta g(s)\right)^2},
\\ 
\frac{+(y_1(s)+(y_1)'_g(s)\, \Delta g(s))\, f(s)}{W_g(y_1,y_2)(s)}=&
\exp_g(\lambda;s)^{-1} \frac{f(s)}{1+\lambda \Delta g(s)}.
\end{aligned}
\end{displaymath}
Hence, Theorem~\ref{solnohomoconog} gives the solution
\begin{equation}\label{solC}
\begin{aligned}
v(t)=&x_0 \, \exp_g(\lambda;t)+\left(v_0-\lambda x_0\right) \varphi(t) \exp_g(\lambda;t)\\
&- \exp_g(\lambda;t) \int_{[0,t)} 
\exp_g(\lambda;s)^{-1}\,
\left[
\varphi(s)\, \frac{f(s)}{1+\lambda \Delta g(s)}+
\frac{f(s)\, \Delta g(s)}{\left(1+\lambda \Delta g(s)\right)^2}
 \right]\, \operatorname{d} \mu_g(s)\\
 &+\varphi(t)\, \exp_g(\lambda;t) \, \int_{[0,t)} \exp_g(\lambda;s)^{-1} \frac{f(s)}{1+\lambda \Delta g(s)}
 \, \operatorname{d} \mu_g(s),
\end{aligned}
\end{equation}
which finishes the result.
\end{enumerate}
\end{proof}

\begin{rem} 
	Observe that the general solution in \eqref{general} coincides 
with the one in Corollary~\ref{constantcoeffsol}. Indeed, we show that this is the case going case by case.
\begin{enumerate}
\item If $\lambda_1\neq \lambda_2$, the solution in \eqref{general} 
can be expressed as
\begin{displaymath}
\begin{aligned}
v(t)=&x_0\exp_g(\lambda_2;t) +
\frac{v_0-\lambda_2x_0}{\lambda_1-\lambda_2} \,\exp_g(\lambda_2;t)\,\left[
 \exp_g
\left( \frac{\lambda_1-\lambda_2}{1+\lambda_2\, \Delta g};t\right) -1\right]\\
&+\frac{1}{\lambda_1-\lambda_2}\exp_g(\lambda_2;t) \int_{[0,t)} 
\exp_g
\left( \frac{\lambda_1-\lambda_2}{1+\lambda_2\, \Delta g};\cdot\right)'_g(s) \\
& \cdot
\bigg(\int_{[0,s)} \exp_g(\lambda_1;r)^{-1}\, \frac{f(r)}{1+\lambda_1 \Delta g(r)}\, \operatorname{d} \mu_g(r)\bigg)\, \operatorname{d} \mu_g(s).
\end{aligned}
\end{displaymath}
Now, if in Lemma~\ref{InByPart} we take the functions given by
\[
w_1(t)=\exp_g
\left( \frac{\lambda_1-\lambda_2}{1+\lambda_2\, \Delta g};t\right),\quad
w_2(t)=\int_{[0,t)} \exp_g(\lambda_1;r)^{-1}\, \frac{f(r)}{1+\lambda_1 \Delta g(r)}\, \operatorname{d} \mu_g(r),
\]
where $t\in[0,T]$, we have, thanks to~\eqref{medtneqts} and Proposition~\ref{gexpprop}, that
\begin{displaymath}
\begin{aligned}
&\int_{[0,t)} 
\exp_g
\left( \frac{\lambda_1-\lambda_2}{1+\lambda_2\, \Delta g};\cdot\right)'_g(s) \cdot
\bigg(\int_{[0,s)} \exp_g(\lambda_1;r)^{-1}\, \frac{f(r)}{1+\lambda_1 \Delta g(r)}\, \operatorname{d} \mu_g(r)\bigg)\, \operatorname{d} \mu_g(s) \\
=&
\exp_g
\left( \frac{\lambda_1-\lambda_2}{1+\lambda_2\, \Delta g};t\right)
\int_{[0,t)} \exp_g(\lambda_1;s)^{-1}\, \frac{f(s)}{1+\lambda_1 \Delta g(s)}\, \operatorname{d} \mu_g(s)\\
&-\int_{[0,t)} \exp_g(\lambda_2;s)^{-1} \frac{f(s)}{1+\lambda_2 \Delta g(s)}
\, \operatorname{d} \mu_g(s).
\end{aligned}
\end{displaymath}
Therefore,
\begin{equation*}\label{solgA}
\begin{aligned}
v(t)=&\frac{v_0-\lambda_2x_0}{\lambda_1-\lambda_2} \,
\exp_g(\lambda_1;t)+\frac{\lambda_1x_0-v_0}{\lambda_1-\lambda_2} \,\exp_g(\lambda_2;t)\\
&+\frac{1}{\lambda_1-\lambda_2} \exp_g(\lambda_1;t)\int_{[0,t)} \exp_g(\lambda_1;s)^{-1}\, \frac{f(s)}{1+\lambda_1 \Delta g(s)}\, \operatorname{d} \mu_g(s)\\
&-\frac{1}{\lambda_1-\lambda_2} \exp_g(\lambda_2;t)\int_{[0,t)} \exp_g(\lambda_2;s)^{-1} \frac{f(s)}{1+\lambda_2 \Delta g(s)}
\, \operatorname{d} \mu_g(s),
\end{aligned}
\end{equation*}
which coincides with \eqref{solgB}.
\item On the other hand, if $\lambda_1=\lambda_2=\lambda$, the expression~\eqref{general} 
can be rewritten as
\begin{displaymath}
\begin{aligned}
v(t)=&x_0\exp_g(\lambda;t) +
\left( v_0-\lambda x_0\right) \,\exp_g(\lambda;t)\, \int_{[0,t)} \frac{1}{1+\lambda \Delta g(s)}\, \operatorname{d} \mu_g(s)\\
&+\exp_g(\lambda;t)\int_{[0,t)} \frac{1}{1+\lambda \Delta g(s)}\, 
\bigg(\int_{[0,s)} \exp_g(\lambda;r)^{-1}\, \frac{f(r)}{1+\lambda \Delta g(r)}\, \operatorname{d} \mu_g(r)\bigg)\, \operatorname{d} \mu_g(s).
\end{aligned}
\end{displaymath}
Once again, applying Lemma~\ref{InByPart} to the functions
\[
w_1(t)= \int_{[0,t)} \frac{1}{1+\lambda \Delta g(s)}\, \operatorname{d}\mu_g(s), \quad w_2(t)=\int_{[0,t)} \exp_g(\lambda;s)^{-1}\, \frac{f(s)}{1+\lambda \Delta g(s)}\, \operatorname{d} \mu_g(s),
\]
where  $t\in[0,T]$, we have, thanks to~\eqref{medtneqts} and Proposition~\ref{gexpprop}, that
\begin{displaymath}
\begin{aligned}
& \int_{[0,t)} \frac{1}{1+\lambda \Delta g(s)}\, 
\bigg(\int_{[0,s)} \exp_g(\lambda;r)^{-1}\, \frac{f(r)}{1+\lambda \Delta g(r)}\, \operatorname{d} \mu_g(r)\bigg)\, \operatorname{d} \mu_g(s)\\
=&\left(\int_{[0,t)} \frac{1}{1+\lambda \Delta g(s)}, \operatorname{d} \mu_g(s)\right)
\bigg(\int_{[0,t)} \exp_g(\lambda;s)^{-1}\, \frac{f(s)}{1+\lambda \Delta g(s)}\, \operatorname{d} \mu_g(s)\bigg)
\\
&-\int_{[0,t)} 
\bigg(\int_{[0,s)} \frac{1}{1+\lambda \Delta g(r)} \, \operatorname{d} \mu_g(r)\bigg)
\exp_g(\lambda;s)^{-1}\, \frac{f(s)}{1+\lambda \Delta g(s)}\, \operatorname{d} \mu_g(s)\\
&-\int_{[0,t)} \exp_g(\lambda;s)^{-1} \frac{f(s)\Delta g(s)}{(1+\lambda \Delta g(s))^2}\, 
\operatorname{d} \mu_g(s).\\
\end{aligned}
\end{displaymath}
Therefore,
\begin{displaymath}
\begin{aligned}
v(t)=&x_0\exp_g(\lambda;t) +
\left( v_0-\lambda x_0\right) \,\exp_g(\lambda;t)\, \int_{[0,t)} \frac{1}{1+\lambda \Delta g(s)}\, \operatorname{d} \mu_g(s)\\
&+\exp_g(\lambda;t)\left(\int_{[0,t)} \frac{1}{1+\lambda \Delta g(s)}\, \operatorname{d} \mu_g(s)\right)
\bigg(\int_{[0,t)} \exp_g(\lambda;s)^{-1}\, \frac{f(s)}{1+\lambda \Delta g(s)}\, \operatorname{d} \mu_g(s)\bigg)
\\
&-\exp_g(\lambda;t) \int_{[0,t)} \Bigg[
\bigg(\int_{[0,s)} \frac{1}{1+\lambda \Delta g(r)} \, \operatorname{d} \mu_g(r)\bigg)
\exp_g(\lambda;s)^{-1}\, \frac{f(s)}{1+\lambda \Delta g(s)}\\
&+\exp_g(\lambda;s)^{-1} \frac{f(s)\Delta g(s)}{(1+\lambda \Delta g(s))^2}\Bigg]\, 
\operatorname{d} \mu_g(s).
\end{aligned}
\end{displaymath}
\end{enumerate}
\end{rem}

Several numerical examples can be found in 
\cite[Section 6]{Fernandez2021}. In this section the authors 
present an application to the Stieltjes harmonic oscillator as well as 
an example in which the resonance effect appears.

\subsection{An example with variable coefficients}

In this section we will see how we can apply the theory we have developed to the following problem,
\begin{empheq}[left=\empheqlbrace]{align} 
& v_g''(t)+P(t)\, v_g'(t)=f(t),\quad  t \in [0,T],
\label{ejenocteA}\\
& v(0)=x_0,\;v_g'(0)=v_0, \label{ejenocteB}
\end{empheq}
where $x_0,\, v_0 \in \mathbb{C}$, $P,\, f\in \mathcal{BC}_g([0,T];\mathbb{C})$ are such that $1-P(t)\, \Delta g(t)\neq 0$, $ t\in [0,T]\cap D_g$. Observe that in this case it is clear that $y_1(t)=1$, $t\in [0,T]$, is a solution of the 
homogeneous equation $v_g''(t)-P(t)\, v_g'(t)=0$, $ t \in [0,T]$. Therefore, we can use Theorem~\ref{solhomoconog} to obtain another linearly independent solution of the homogeneous equation, 
$y_2=v\, y_1$, with
\begin{displaymath}
v(t)=\int_{[0,t)} \exp_g(-P;s)\, \operatorname{d}\mu_g(s).
\end{displaymath}
Thus, the solution of the initial value problem~\eqref{ejenocteA}-\eqref{ejenocteB} can be obtained thanks to Theorem~\ref{solnohomoconog}, which yields
\begin{displaymath}
\begin{aligned}
v(t)
=& \left(\frac{x_0\, (y_2)'_g(0)-v_0\, y_2(0)}{W_g(y_1,y_2)(0)}\right) \, y_1(t)+
\left(\frac{v_0\, y_1(0)-x_0\, (y_1)'_g(0)}{W_g(y_1,y_2)(0)}\right) \, y_2(t)\\
&+
y_1(t)\, \int_{[0,t)} \frac{-(y_2+(y_2)'_g\, \Delta g)\, f}{W_g(y_1,y_2)}\, \operatorname{d}\mu_g+y_2(t)\, \int_{[0,t)} \frac{(y_1+(y_1)'_g\, \Delta g)\, f}{W_g(y_1,y_2)}\, \operatorname{d}\mu_g.
\end{aligned}
\end{displaymath}
In this case, we have that for $t\in[0,T]$,
\begin{displaymath}
\begin{aligned}
W_g(y_1,y_2)(t)=&y_1(t) \, (y_2)'_g(t) + y_1(t)\, (y_2)''_g(t) \, \Delta g(t)\\
=& \exp_g(-P;t)-P(t)\, \exp_g(-P;t)\, \Delta g(t)= \exp_g(-P;t)\left(1-P(t)\, \Delta g(t)\right)\neq 0,
\end{aligned}
\end{displaymath}
and, as a consequence,
\begin{align*}
-\frac{y_2(t)+(y_2)'_g(t)\, \Delta g(t)}{W_g(y_1,y_2)(t)}
&=-\frac{
\int_{[0,t)} \exp_g(-P;s)\, \operatorname{d}\mu_g(s)+\exp_g(-P;t)\, \Delta g(t)
}{ \exp_g(-P;t)\left(1-P(t)\, \Delta g(t)\right)},\\
 \frac{y_1(t)+(y_1)'_g(t)\, \Delta g(t)}{W_g(y_1,y_2)(t)}&
 =\frac{1}{\exp_g(-P;t)\left(1-P(t)\, \Delta g(t)\right)}.
\end{align*}
Therefore,
\begin{displaymath}
\begin{aligned}
v(t)=&x_0+v_0 \, \int_{[0,t)} \exp_g(-P;s)\, \operatorname{d}\mu_g(s)\\
&-\int_{[0,t)}\left( \int_{[0,s)} \exp_g(-P;u)\,\operatorname{d}\mu_g(u)\right)\,\exp_g(-P;s)^{-1} \frac{f(s)}{1-P(s)\, \Delta g(s)}\, \operatorname{d}\mu_g(s)\\
&-\int_{[0,t)}  \frac{f(s)\, \Delta g(s)}{1-P(s)\, \Delta g(s)}\, \operatorname{d}\mu_g(s)\\
&+\left(\int_{[0,t)} \exp_g(-P;s)\, \operatorname{d}\mu_g(s)\right)\, \left(\int_{[0,t)} 
\exp_g(-P;s)^{-1} \frac{f(s)}{1-P(s)\, \Delta g(s)}\, \operatorname{d}\mu_g(s)\right).
\end{aligned}
\end{displaymath}
By using integration by parts (Lemma~\ref{InByPart}),
\begin{displaymath}
\begin{aligned}
&-\int_{[0,t)}\left( \int_{[0,s)} \exp_g(-P;u)\,\operatorname{d}\mu_g(u)\right)\,\exp_g(-P;s)^{-1} \frac{f(s)}{1-P(s)\, \Delta g(s)}\, \operatorname{d}\mu_g(s)\\
=&\int_{[0,t)} \exp_g(P;s) \left(\int_{[0,s)} \exp_g(-P;u)^{-1}\, \frac{f(u)}{1-P(u)\, \Delta g(u)} \right)\, \operatorname{d}\mu_g(s)\\
&+\int_{[0,t)} \exp_g(P;s)\exp_g(-P;s)^{-1} \frac{f(s)}{1-P(s)\, \Delta g(s)}\, \Delta g(s)\, \operatorname{d}\mu_g(s)\\
&-\left(\int_{[0,t)}  \exp_g(P;s)\, \operatorname{d}\mu_g(s) \right)\left(\int_{[0,t)}  \exp_g(-P;s)^{-1} \frac{f(s)}{1-P(s)\, \Delta g(s)}\, \operatorname{d}\mu_g(s) \right).
\end{aligned}
\end{displaymath}
Thus,
\begin{displaymath}
\begin{aligned}
v(t)=&x_0+v_0 \, \int_{[0,t)} \exp_g(-P;s)\, \operatorname{d}\mu_g(s)\\
&+\int_{[0,t)} \exp_g(P;s) \left(\int_{[0,s)} \exp_g(-P;u)^{-1}\, \frac{f(u)}{1-P(u)\, \Delta g(u)} \right)\, \operatorname{d}\mu_g(s).
\end{aligned}
\end{displaymath}
Note that the previous formula is the same as the one we obtain if we introduce the variable $w=v_g'$ and 
use~\cite[Proposition 4.12]{Fernandez2021}. Indeed, we have that $w$ satisfies the initial value problem
\begin{displaymath}
\left\{\begin{aligned}
& w_g'(t)+P(t)\, w(t)=f(t)\quad t\in [0,T],\\
& w(0)=v_0.
\end{aligned}\right.
\end{displaymath}
Therefore, by \cite[Proposition 4.12]{Fernandez2021},
\begin{displaymath}
w(t)=v_0\, \exp_g(-P;t)+\exp_g(-P;t)\, \int_{[0,t)} \exp_g(-P;s)^{-1} \frac{f(s)}{1-P(s)\, \Delta g(s)}\,\operatorname{d}\mu_g(s),
\end{displaymath}
and, as a consequence,
\begin{displaymath}
\begin{aligned}
v(t)=&x_0+\int_{[0,t)} w(s)\, \operatorname{d}\mu_g(s)\\
=& x_0 + v_0\, \int_{[0,t)}  \exp_g(-P;s)\, \operatorname{d}\mu_g(s)\\
&+\int_{[0,t)} \exp_g(P;s) \left(\int_{[0,s)} \exp_g(-P;u)^{-1}\, \frac{f(u)}{1-P(u)\, \Delta g(u)} \right)\, \operatorname{d}\mu_g(s).
\end{aligned}
\end{displaymath}

\subsection{One-dimensional linear Helmholtz equation with 
piecewise-constant coefficients}

In the two previous examples we have seen how the theory 
we have developed allows us to recover the same solution 
as the one obtained using order reduction methods. In 
the following example we will see that it is possible 
to use the techniques that we have developed throughout this paper to study the one-dimensional linear Helmholtz equation 
with piecewise-constant coefficients given by
\begin{empheq}[left=\empheqlbrace]{align} 
& v_g''(t)+w_0(t)^2\, v(t)=f(t),\quad  t \in [0,T],
\label{eq:piecewisenhA}\\
& v(0)=x_0,\;v_g'(0)=v_0, \label{eq:piecewisenhB}
\end{empheq}
where $x_0,\, v_0 \in \mathbb{C}$, $f\in \mathcal{BC}_g([0,T];\mathbb{C})$, 
and 
\begin{displaymath}
w_0:t\in [0,T]\rightarrow w_0(t)=\left\{
\begin{aligned}
&w_{1}, & t\in [0,t_1], \\
&w_{2}, & t\in (t_1,T],
\end{aligned} \right.
\end{displaymath}
for some $t_1\in [0,T]\cap D_g$ and $w_1,\, w_2\in \mathbb{R}$. This equation, in the case of the usual derivative, was studied in \cite{Bugarija2020} in the context of inverse scattering problems. Observe that
\begin{equation}\label{gcond0}
\left(1\pm i\, w_1\, \Delta g(t)\right)\,\left(1\pm i\, w_2\, 
\Delta g(t)\right)\neq 0,\quad t \in [0,T]\cap D_g.
\end{equation}
We have that $w_0\in \mathcal{BC}_g([0,T];\mathbb{C})$ and 
$\exp_g(\pm i\, w_k;\cdot)\in \mathcal{BC}_g^{\infty}([0,T];
\mathbb{C})$, $k=1,2$, are well-defined and nonzero on $[0,T]$. 

We have, thanks 
to Remark~\ref{remexistence}, 
that~\eqref{eq:piecewisenhA}-\eqref{eq:piecewisenhB} has an 
unique solution $v\in \mathcal{BC}_g^2([0,T];\mathbb{C})$. Let us 
see now that we can find an explicit solution of 
the problem~\eqref{eq:piecewisenhA}-\eqref{eq:piecewisenhB} using 
the method of variation of parameters. Let us then consider the homogeneous 
problem
\begin{empheq}[left=\empheqlbrace]{align} 
& v_g''(t)+w_0(t)^2\, v(t)=0,\quad  t \in [0,T],
\label{eq:piecewisehA}\\
& v(0)=x_0,\;v_g'(0)=v_0, \label{eq:piecewisehB}
\end{empheq}
and denote by
\[\lambda_1^1:=i\,w_1,\quad \lambda_1^2:=-i\, w_1,\quad \lambda_2^1:= i\,w_2,\quad \lambda_2^2:=-i\, w_2.
\]
Let us see that it is possible to obtain two linearly 
independent solutions of \eqref{eq:piecewisehA}
of the form 
\begin{displaymath}
y_k(t)=\exp_g(\lambda_1^k;t)\, \chi_{[0,t_1]}(t)+
\left(\alpha_1^k\, \exp_g(\lambda_2^1;t)+\alpha_2^k\,
\exp_g(\lambda_2^2;t) \right)\,\chi_{(t_1,T]}(t),\quad t\in [0,T],\quad k=1,2,
\end{displaymath}
for some $\alpha_1^k,\,\alpha_2^k\in \mathbb{C}$, $k=1,2$. We have that 
$\chi_{[0,t_1]},\,\chi_{(t_1,T]}\in \mathcal{BC}_g([0,T];\mathbb{C})$. Thus, $y_k\in \mathcal{BC}_g([0,T];\mathbb{C})$. Moreover, thanks to 
Proposition \ref{PropStiDer} and the properties of the $g$-exponential we have that 
$(y_k)''_g(t)+w_0(t)^2\, y_k(t)=0$, $t\in [0,T]\setminus \{t_1\}$, 
$k=1,2$. Now, for $k=1,2,$ bearing in mind that $\exp_g(p;t^+)=(1+p(t^+)\Delta g(t))\exp_g(p;t^+)$ for any $t\in D_g$, c.f. \cite[Remark 3.3]{Ma21}, we have that
	\begin{align*}
		(y_k)'_g(t_1)&=\frac{y_k(t_1^+)-y_k(t_1)}{\Delta g(t_1)}\\
		&=\frac{\alpha_1^k(1+\lambda_2^1\Delta g(t_1)) \exp_g(\lambda_2^1;t_1)+\alpha_2^k(1+\lambda_2^2\Delta g(t_1))
			\exp_g(\lambda_2^2;t_1)
			-\exp_g(\lambda_1^k;t_1)}{\Delta g(t_1)}.
	\end{align*}
Therefore, if we want to ensure that the first $g$-derivative is $g$-continuous 
at $t=t_1$, it must be satisfied that $\lim_{t\to t_1^-} (y_k)'_g(t)=(y_k)'_g(t_1)$, $k=1,2,$ see Proposition~\ref{proreg}.
In particular, since the $g$-exponential is $g$-continuous, it is left-continuous and thus, for $k=1,2$, we require that
\begin{displaymath}
\begin{aligned}
\lambda_1^k\,\exp_g(\lambda_1^k;t_1)=&
\frac{\alpha_1^k(1+\lambda_2^1\Delta g(t_1)) \exp_g(\lambda_2^1;t_1)+\alpha_2^k(1+\lambda_2^2\Delta g(t_1))
	\exp_g(\lambda_2^2;t_1)
	-\exp_g(\lambda_1^k;t_1)}{\Delta g(t_1)}
\end{aligned}
\end{displaymath}
or, equivalently,
\begin{equation}\label{gcond1}
\left(1+\lambda_1^k \Delta g(t_1)\right) \exp_g(\lambda_1^k;t_1)=
\alpha_1^k\exp_g(\lambda_2^1;t_1)\left(1+\lambda_2^1\Delta g(t_1) \right)+\alpha_2^k\exp_g(\lambda_2^2;t_1)\left(1+\lambda_2^2\,\Delta g(t_1) \right).
\end{equation}
In the case in which~\eqref{gcond1} is fulfilled, we have that for $k=1,2$ and $t\in [0,T]$,
\begin{displaymath}
(y_k)_g'(t)=\lambda_1^k\,\exp_g(\lambda_1^k;t)\, \chi_{[0,t_1]}(t)+
\left(\alpha_1^k\,\lambda_2^1\, \exp_g(\lambda_2^1;t)+\alpha_2^k\,
\lambda_2^2\,
\exp_g(\lambda_2^2;t) \right)\,\chi_{(t_1,T]}(t).
\end{displaymath}
Observe that $(y_k)_g'\in \mathcal{BC}_g([0,T];\mathbb{C})$.

In a similar fashion, for $k=1,2$, we compute the
second $g$-derivative at the point $t=t_1$,
\begin{displaymath}
\begin{aligned}
(y_k)_g''(t_1)=&\frac{(y_k)'_g(t_1^+)-(y_k)'_g(t_1)}{\Delta g(t_1)}\\
=& \frac{
\alpha_1^k(1+\lambda_2^1\Delta g(t_1))\,\lambda_2^1\, \exp_g(\lambda_2^1;t_1)+
\alpha_2^k(1+\lambda_2^2\Delta g(t_1))\, \lambda_2^2\,
\exp_g(\lambda_2^2;t_1)
-\lambda_1^k\,\exp_g(\lambda_1^k;t_1)}{\Delta g(t_1)}.
\end{aligned}
\end{displaymath}
Reasoning analogously to the previous case, for $k=1,2$, in order for the second $g$-derivative of $y_k$ to be $g$-continuous at the point $t=t_1$, it must be satisfied 
that
\begin{displaymath}
\begin{aligned}
(\lambda_1^k)^2\, \exp_g(\lambda_1^k;t_1)=&
\frac{1}{\Delta g(t_1)} \bigg(-\lambda_1^k\,\exp_g(\lambda_1^k;t_1) 
+\alpha_1^k\,\lambda_2^1\, \exp_g(\lambda_2^1;t_1)\,\left(
1+\lambda_2^1\, \Delta g(t_1)
\right)\\&+\alpha_2^k\,\lambda_2^2\, \exp_g(\lambda_2^2;t_1)\, 
\left( 1+\lambda_2^2\,\Delta g(t_1)\right)
\bigg).
\end{aligned}
\end{displaymath} 
Again, for $k=1,2$,  this is equivalent to
\begin{equation}\label{gcond2}
\begin{aligned}
\left(1+\lambda_1^k\, \Delta g(t_1)\right)\,\lambda_1^k\, \exp_g(\lambda_1^k;t_1)=&
\alpha_1^k\,\lambda_2^1\,\exp_g(\lambda_2^1;t_1)\,\left(1+\lambda_2^1\,\Delta g(t_1) \right)\\&+\alpha_2^k\,\lambda_2^2\,\exp_g(\lambda_2^2;t_1)\,\left(1+\lambda_2^2\,\Delta g(t_1) \right),
\end{aligned}
\end{equation}
so, in the case in which~\eqref{gcond2} is fulfilled, we obtain, for $t\in [0,T]$,
\begin{displaymath}
\begin{aligned}
(y_k)_g''(t)=&(\lambda_1^k)^2\,\exp_g(\lambda_1^k;t)\, \chi_{[0,t_1]}(t)+
\left(\alpha_1^k\,(\lambda_2^1)^2\, \exp_g(\lambda_2^1;t)+\alpha_2^k\,
(\lambda_2^2)^2\,
\exp_g(\lambda_2^2;t) \right)\,\chi_{(t_1,T]}(t).
\end{aligned}
\end{displaymath}
Observe that $(y_k)_g''\in \mathcal{BC}_g([0,T];\mathbb{C})$.

Conditions~\eqref{gcond1} and~\eqref{gcond2} give us 
the following linear system  of equations
\begin{equation} \label{eq:systemalphak}
\begin{aligned}
&\left(\begin{array}{rr}
\exp_g(\lambda_2^1;t_1)\,\left(1+\lambda_2^1\,\Delta g(t_1) \right)&
\exp_g(\lambda_2^2;t_1)\,\left(1+\lambda_2^2\,\Delta g(t_1) \right)\\
\lambda_2^1\, \exp_g(\lambda_2^1;t_1)\,\left(
1+\lambda_2^1\, \Delta g(t_1)
\right)& 
\lambda_2^2\,\exp_g(\lambda_2^2;t_1)\,\left(1+\lambda_2^2\,\Delta g(t_1) \right)
\end{array}\right)
\left(\begin{array}{c}
\alpha_1^k\\\alpha_2^k
\end{array}
\right)\\
=&
\left(\begin{array}{r}
\left(1+\lambda_1^k\, \Delta g(t_1)\right)\, \exp_g(\lambda_1^k;t_1)\\
\lambda_1^k\,\left(1+\lambda_1^k\, \Delta g(t_1)\right)\, \exp_g(\lambda_1^k;t_1)
\end{array}
\right).
\end{aligned}
\end{equation}
Note that, thanks to~\eqref{gcond0}, the determinant of the matrix $A\in \mathcal{M}_{2\times 2}(
\mathbb{C})$ of the previous system is such that
\begin{displaymath}
\det(A)=\exp_g(\lambda_2^1;t_1)
\, \exp_g(\lambda_2^2;t_1)\,\left(1+\lambda_2^1\,\Delta g(t_1) \right)\,\left(1+\lambda_2^2\,\Delta g(t_1) \right)\,\left(
\lambda_2^2-\lambda_2^1\right)\neq 0.
\end{displaymath}
Therefore,
\begin{equation}\label{defalphak1}
	\begin{aligned}
	\alpha_1^k =&\frac{\exp_g(\lambda_1^k;t_1)\,\exp_g(\lambda_2^2;t_1)\,\left(1+\lambda_1^k\, \Delta g(t_1)\right)\,\left(1+\lambda_2^2\,\Delta g(t_1) \right)\,(\lambda_2^2-\lambda_1^k)}{\exp_g(\lambda_2^1;t_1)
		\, \exp_g(\lambda_2^2;t_1)\,\left(1+\lambda_2^1\,\Delta g(t_1) \right)\,\left(1+\lambda_2^2\,\Delta g(t_1) \right)\,\left(
		\lambda_2^2-\lambda_2^1\right)}\\
	=&\frac{\exp_g(\lambda_1^k;t_1)\,\left(1+\lambda_1^k\, \Delta g(t_1)\right)\,(\lambda_2^2-\lambda_1^k)}{\exp_g(\lambda_2^1;t_1)
		\,\left(1+\lambda_2^1\,\Delta g(t_1)\right)\, \left(
		\lambda_2^2-\lambda_2^1\right)},\end{aligned}
	\end{equation}
	\begin{equation}\label{defalphak2}
		\begin{aligned}
	\alpha_2^k=&\frac{\exp_g(\lambda_2^1;t_1)\,\exp_g(\lambda_1^k;t_1)\,\left(1+\lambda_2^1\,\Delta g(t_1) \right)\,\left(1+\lambda_1^k\, \Delta g(t_1)\right)\, (\lambda_1^k-\lambda_2^1)}{\exp_g(\lambda_2^1;t_1)
		\, \exp_g(\lambda_2^2;t_1)\,\left(1+\lambda_2^1\,\Delta g(t_1) \right)\,\left(1+\lambda_2^2\,\Delta g(t_1) \right)\,\left(
		\lambda_2^2-\lambda_2^1\right)}\\
	=&\frac{\exp_g(\lambda_1^k;t_1)\,\left(1+\lambda_1^k\, \Delta g(t_1)\right)\,(\lambda_1^k-\lambda_2^1)}{\exp_g(\lambda_2^2;t_1)\,\left(1+\lambda_2^2\,\Delta g(t_1) \right)\,\left(
		\lambda_2^2-\lambda_2^1\right)}.
	\end{aligned}\end{equation}
	
We have the following result.

\begin{cor}\label{corollaryyk}
	Let $t_1\in [0,T]\cap D_g$, $w_1,\,w_2\in 
\mathbb{R}$ and
consider $\alpha_1^k,\,\alpha_2^k\in \mathbb{C}$ as in \eqref{defalphak1}-\eqref{defalphak2}. Then, the maps 
\begin{displaymath}
y_k=\exp_g(\lambda_1^k;\cdot)\, \chi_{[0,t_1]}+
\left(\alpha_1^k\, \exp_g(\lambda_2^1;\cdot)+\alpha_2^k\,
\exp_g(\lambda_2^2;\cdot) \right)\,\chi_{(t_1,T]},\quad k=1,2,
\end{displaymath}
are two linearly independent solutions of the homogeneous equation~\eqref{eq:piecewisehA}. In particular, 
\begin{equation}\label{eq:solpiecewiseh}
\begin{aligned}
v_h(t)=&
\left(\frac{\lambda_1^2\, x_0-v_0}{\lambda_1^2-\lambda_1^1} \right)\,y_1(t)
+\left(\frac{v_0-\lambda_1^1\,x_0}{\lambda_1^2-\lambda_1^1} \right)\,y_2(t),
\end{aligned}
\end{equation}
is the solution of the initial value problem~\eqref{eq:piecewisehA}-\eqref{eq:piecewisehB}. 
\end{cor}

\begin{proof} 
Given the reasoning above, it is clear that $y_1,\, y_2$ are two solutions of \eqref{eq:piecewisehA}, so let us show that
$y_1,\, y_2$ are linearly independent. Indeed, by Lemma~\ref{GwronRel}, we have that
\begin{displaymath}
W_g(y_1,y_2)=\left(1+w_0^2\, \Delta g^2\right) 
\widetilde{W}_g(y_1,y_2)=\left(1+w_0^2\, \Delta g^2\right) 
\widetilde{W}_g(y_1,y_2)(0)\,
\exp_g(w_0^2\,\Delta g;\cdot).
\end{displaymath}
 Now, given that by definition $\widetilde{W}_g(y_1,y_2)=y_1(y_2)'_g - y_2 (y_1)'_g$,
we have, by condition \eqref{gcond0}, that
\begin{displaymath}
W_g(y_1,y_2)=(\lambda_1^2-\lambda_1^1)\,\left(1+w_0^2\, \Delta g^2\right)\,
\exp_g(w_0^2\,\Delta g;\cdot)\neq 0,
\end{displaymath}
so Lemma~\ref{lemwr} ensures that $y_1,y_2$ are linearly independent.
\end{proof}

We have the following result for the non homogeneous 
problem~\eqref{eq:piecewisenhA}-\eqref{eq:piecewisenhB}.

\begin{cor} Let $f\in \mathcal{BC}_g([0,T];\mathbb{C})$, 
$t_1\in [0,T]\cap D_g$, $w_1,\,w_2\in 
\mathbb{R}$,  and
consider $\alpha_1^k,\,\alpha_2^k\in \mathbb{C}$ as in \eqref{defalphak1}-\eqref{defalphak2}. Then, the solution of 
the initial value problem~\eqref{eq:piecewisenhA}-\eqref{eq:piecewisenhB} is given by
\begin{displaymath}
v=v_p+v_h,
\end{displaymath}
where $v_h$ is the solution of the homogeneous 
problem~\eqref{eq:piecewisehA}-\eqref{eq:piecewisehB} 
given by~\eqref{eq:solpiecewiseh}, and 
\begin{displaymath}
v_p=v_{p_1}\,\chi_{[0,t_1]}+v_{p_2}\, \chi_{(t_1,T]},
\end{displaymath}
is a particular solution of~\eqref{eq:piecewisenhA} with
\begin{displaymath}
\begin{aligned}
v_{p_1}(t)=&\frac{y_1(t)}{(\lambda_1^2-\lambda_1^1)}\, \int_{[0,t)} \frac{-(y_2+(y_2)'_g\, \Delta g)\, f}{\left(1+w_0^2\, \Delta g^2\right)\,
\exp_g(w_0^2\,\Delta g;\cdot)}\, \operatorname{d}\mu_g\\
& +\frac{y_2(t)}{(\lambda_1^2-\lambda_1^1)}\, \int_{[0,t)} \frac{(y_1+(y_1)'_g\, \Delta g)\, f}{\left(1+w_0^2\, \Delta g^2\right)\,
\exp_g(w_0^2\,\Delta g;\cdot)}\, \operatorname{d}\mu_g,\quad t\in[0,T],
\end{aligned}
\end{displaymath}
and, for each $t\in[0,T]$, 
\begin{displaymath}
\begin{aligned}
v_{p_2}(t)=&\frac{\left(\alpha_1^1\, \exp_g(\lambda_2^1;t)
+\alpha_2^1\,\exp_g(\lambda_2^2;t) \right)}{\lambda_1^1-\lambda_1^2} 
\Bigg[
\int_{[0,t_1)} \exp_g(\lambda_1^1;\cdot)^{-1}\, 
\frac{f}{1+\lambda_1^1 \Delta g}\, \operatorname{d} \mu_g\\
&+\exp_g(\lambda_1^1;t_1)^{-1}\, 
\frac{f(t_1)\, \Delta g(t_1)}{1+\lambda_1^1 \Delta g(t_1)}+
\alpha_1^2\,\int_{(t_1,t)} \frac{ \exp_g(\lambda_2^1;\cdot)}{\left(1+\lambda_2^2\, \Delta g\right)\,\exp_g(w_0^2\,\Delta g;\cdot)}\, f\, \operatorname{d}\mu_g\\
&+\alpha_2^2\,
\int_{(t_1,t)}
 \frac{\exp_g(\lambda_2^2;\cdot)}{\left(1+\lambda_2^1\, \Delta g\right)\,\exp_g(w_0^2\,\Delta g;\cdot)}\, f\, \operatorname{d}\mu_g
\Bigg]\\
&-\frac{\left(\alpha_1^2\, \exp_g(\lambda_2^1;t)
+\alpha_2^2\,\exp_g(\lambda_2^2;t) \right)}{\lambda_1^1-\lambda_1^2}
\Bigg[
\int_{[0,t_1)} \exp_g(\lambda_1^2;\cdot)^{-1} 
\frac{f}{1+\lambda_1^2 \Delta g}
\, \operatorname{d} \mu_g\\
&+\exp_g(\lambda_1^2;t_1)^{-1} 
\frac{f(t_1)\, \Delta g(t_1)}{1+\lambda_1^2 \Delta g(t_1)}
+\alpha_1^1\,\int_{(t_1,t)}\frac{ \exp_g(\lambda_2^1;\cdot)}{\left(1+\lambda_2^2\, \Delta g\right)\,\exp_g(w_0^2\,\Delta g;\cdot)}\, f\, \operatorname{d}\mu_g\\
&+
\alpha_2^1\,\int_{(t_1,t)} \frac{\exp_g(\lambda_2^2;\cdot)}{\left(1+\lambda_2^1\, \Delta g\right)\,\exp_g(w_0^2\,\Delta g;\cdot)}\, f\, \operatorname{d}\mu_g\Bigg].
\end{aligned}
\end{displaymath}

\end{cor}

\begin{proof} First, observe that the hypotheses of Theorem~\ref{solnohomoconog} are satisfied for $y_1$ and $y_2$ given by Corollary~\ref{corollaryyk}. Hence, we have that
\begin{displaymath}
\begin{aligned}
v_p(t)=& \frac{y_1(t)}{\lambda_1^2-\lambda_1^1}\, \int_{[0,t)} \frac{-(y_2+(y_2)'_g\, \Delta g)\, f}{\left(1+w_0^2\, \Delta g^2\right)\,
\exp_g(w_0^2\,\Delta g;\cdot)}\, \operatorname{d}\mu_g\\
& +\frac{y_2(t)}{\lambda_1^2-\lambda_1^1}\, \int_{[0,t)} \frac{(y_1+(y_1)'_g\, \Delta g)\, f}{\left(1+w_0^2\, \Delta g^2\right)\,
\exp_g(w_0^2\,\Delta g;\cdot)}\, \operatorname{d}\mu_g.
\end{aligned}
\end{displaymath}
is a particular solution of~\eqref{eq:piecewisenhA}. Now, observe that, for $t\in [0,t_1]$, 
\begin{displaymath}
\begin{aligned}
& \int_{[0,t)} \frac{-(y_2+(y_2)'_g\, \Delta g)\, f}{\left(1+w_0^2\, \Delta g^2\right)\,
\exp_g(w_0^2\,\Delta g;\cdot)}\, \operatorname{d}\mu_g\\
=&\int_{[0,t)}\frac{-\exp_g(\lambda_1^2;\cdot)\,\left(1+
\lambda_1^2\, \Delta g \right)\, f}{
\left(1+\lambda_1^1\,\lambda_1^2\, \Delta^2g\right)\,
\exp_g(\lambda_1^1\,\lambda_1^2\,\Delta g;\cdot)}\, \operatorname{d}\mu_g\\
=& \int_{[0,t)}\frac{-\exp_g(\lambda_1^2;\cdot)\, f}{
\left(1+\lambda_1^1\, \Delta g\right)\, 
\exp_g(\lambda_1^1\,\lambda_1^2\,\Delta g;\cdot)}\, \operatorname{d}\mu_g\\
=& \int_{[0,t)} -\frac{f}{1+\lambda_1^1\, \Delta g}\,
\exp_g(\lambda_1^2;\cdot)\, \exp_g\left(-\frac{\lambda_1^1\,\lambda_1^2\,
\Delta g}{1+\lambda_1^1\,\lambda_1^2\,\Delta^2 g};\cdot\right)\, 
\operatorname{d}\mu_g\\
=&\int_{[0,t)}-\frac{f}{1+\lambda_1^1\, \Delta g}\,
\exp_g\left(\lambda_1^2-\frac{\lambda_1^1\,\lambda_1^2\,
\Delta g}{1+\lambda_1^1\,\lambda_1^2\,\Delta^2 g} 
-\frac{\lambda_1^1\,(\lambda_1^2)^2\,
\Delta^2 g}{1+\lambda_1^1\,\lambda_1^2\,\Delta^2 g}
;\cdot\right)\, \operatorname{d}\mu_g\\
=&\int_{[0,t)} -\frac{f}{1+\lambda_1^1\, \Delta g}\,
\exp_g\left(\frac{\lambda_1^2-\lambda_1^1\,\lambda_1^2\,
\Delta g}{1+\lambda_1^1\,\lambda_1^2\,\Delta^2 g} 
;\cdot\right)\, \operatorname{d}\mu_g\\
=&\int_{[0,t)} -\frac{f}{1+\lambda_1^1\, \Delta g}\,
\exp_g\left(\frac{-\lambda_1^1+\lambda_1^1\,\lambda_1^1\,
\Delta g}{1-\lambda_1^1\,\lambda_1^1\,\Delta^2 g} 
;\cdot\right)\, \operatorname{d}\mu_g\\
=&\int_{[0,t)}
-\frac{f}{1+\lambda_1^1\, \Delta g}\,\exp_g\left( 
\frac{-\lambda_1^1}{1+\lambda_1^1\, \Delta g};\cdot
\right)\, \operatorname{d}\mu_g\\
=&\int_{[0,t)} -\frac{f}{1+\lambda_1^1\, \Delta g}\,
\exp_g\left(\lambda_1^1;\cdot\right)^{-1}\, \operatorname{d}\mu_g.
\end{aligned}
\end{displaymath}Analogously, for $t\in [0,t_1]$,
\[\int_{[0,t)} \frac{(y_1+(y_1)'_g\, \Delta g)\, f}{\left(1+w_0^2\, \Delta g^2\right)\,\exp_g(w_0^2\,\Delta g;\cdot)}\, \operatorname{d}\mu_g=\int_{[0,t)} \exp_g(\lambda_1^2;\cdot)^{-1} 
\frac{f}{1+\lambda_1^2 \Delta g}
\, \operatorname{d} \mu_g.\]
Therefore,
\begin{displaymath}
\begin{aligned}
v_p|_{[0,t_1]}(t)=&\frac{1}{\lambda_1^1-\lambda_1^2} \exp_g(\lambda_1^1;t)
\int_{[0,t)} \exp_g(\lambda_1^1;\cdot)^{-1}\, 
\frac{f}{1+\lambda_1^1 \Delta g}\, \operatorname{d} \mu_g\\
&-\frac{1}{\lambda_1^1-\lambda_1^2} \exp_g(\lambda_1^2;t)
\int_{[0,t)} \exp_g(\lambda_1^2;\cdot)^{-1} 
\frac{f}{1+\lambda_1^2 \Delta g}
\, \operatorname{d} \mu_g.
\end{aligned}
\end{displaymath}
On the other hand, given $t\in (t_1,T)$,
\begin{displaymath}
\begin{aligned}
v_p(t)=&\frac{y_1(t)}{\lambda_1^1-\lambda_1^2} 
\bigg(
\int_{[0,t_1)} \exp_g(\lambda_1^1;\cdot)^{-1}\, 
\frac{f}{1+\lambda_1^1 \Delta g}\, \operatorname{d} \mu_g+\exp_g(\lambda_1^1;t_1)^{-1}\, 
\frac{f(t_1)\, \Delta g(t_1)}{1+\lambda_1^1 \Delta g(t_1)}\\
&+
\int_{(t_1,t)} \frac{(y_2+(y_2)'_g\, \Delta g)\, f}{\left(1+w_0^2\, \Delta g^2\right)\,
\exp_g(w_0^2\,\Delta g;\cdot)}\, \operatorname{d}\mu_g
\bigg)\\
&-\frac{y_2(t)}{\lambda_1^1-\lambda_1^2}
\bigg(
\int_{[0,t_1)} \exp_g(\lambda_1^2;\cdot)^{-1} 
\frac{f}{1+\lambda_1^2 \Delta g}
\, \operatorname{d} \mu_g+\exp_g(\lambda_1^2;t_1)^{-1} 
\frac{f(t_1)\, \Delta g(t_1)}{1+\lambda_1^2 \Delta g(t_1)}\\
&
+\int_{(t_1,t)}\frac{(y_1+(y_1)'_g\, \Delta g)\, f}{\left(1+w_0^2\, \Delta g^2\right)\,
\exp_g(w_0^2\,\Delta g;\cdot)}\, \operatorname{d}\mu_g\bigg).
\end{aligned}
\end{displaymath}
We have that,
\begin{displaymath}
\begin{aligned}
& \frac{(y_2+(y_2)'_g\, \Delta g)\, f}{\left(1+w_0^2\, \Delta g^2\right)\,
\exp_g(w_0^2\,\Delta g;\cdot)}\\
=&\frac{\left( \alpha_1^2\, \exp_g(\lambda_2^1;\cdot)+\alpha_2^2\,
\exp_g(\lambda_2^2;\cdot)\right)+
\left( \alpha_1^2\,\lambda_2^1\, \exp_g(\lambda_2^1;\cdot)+\alpha_2^2\,
\lambda_2^2\,
\exp_g(\lambda_2^2;\cdot)\right)\,\Delta g}{\left(1+w_0^2\, \Delta g^2\right)\,
\exp_g(w_0^2\,\Delta g;\cdot)}\, f\\
=&\frac{\alpha_1^2\, \exp_g(\lambda_2^1;\cdot)\,\left(1+\lambda_2^1\, \Delta g\right)}{\left(1+w_0^2\, \Delta g^2\right)\,
\exp_g(w_0^2\,\Delta g;\cdot)}\, f+
\frac{\alpha_2^2\,\exp_g(\lambda_2^2;\cdot)\,\left(1+\lambda_2^2\, \Delta g\right)}{\left(1+w_0^2\, \Delta g^2\right)\,
\exp_g(w_0^2\,\Delta g;\cdot)}\, f\\
=&\frac{\alpha_1^2\, \exp_g(\lambda_2^1;\cdot)}{\left(1+\lambda_2^2\, \Delta g\right)\,\exp_g(w_0^2\,\Delta g;\cdot)}\, f+
\frac{\alpha_2^2\,\exp_g(\lambda_2^2;\cdot)}{\left(1+\lambda_2^1\, \Delta g\right)\,\exp_g(w_0^2\,\Delta g;\cdot)}\, f.
\end{aligned}
\end{displaymath}
Analogously,
\begin{displaymath}
\begin{aligned}
& \frac{(y_1+(y_1)'_g\, \Delta g)\, f}{\left(1+w_0^2\, \Delta g^2\right)\,
\exp_g(w_0^2\,\Delta g;\cdot)}=
\frac{\alpha_1^1\, \exp_g(\lambda_2^1;\cdot)}{\left(1+\lambda_2^2\, \Delta g\right)\,\exp_g(w_0^2\,\Delta g;\cdot)}\, f+
\frac{\alpha_2^1\,\exp_g(\lambda_2^2;\cdot)}{\left(1+\lambda_2^1\, \Delta g\right)\,\exp_g(w_0^2\,\Delta g;\cdot)}\, f.
\end{aligned}
\end{displaymath}
Therefore, 
\begin{displaymath}
\begin{aligned}
v_p|_{(t_1,T]}(t)=&\frac{\left(\alpha_1^1\, \exp_g(\lambda_2^1;t)
+\alpha_2^1\,\exp_g(\lambda_2^2;t) \right)}{\lambda_1^1-\lambda_1^2} 
\Bigg[
\int_{[0,t_1)} \exp_g(\lambda_1^1;\cdot)^{-1}\, 
\frac{f}{1+\lambda_1^1 \Delta g}\, \operatorname{d} \mu_g\\
&+\exp_g(\lambda_1^1;t_1)^{-1}\, 
\frac{f(t_1)\, \Delta g(t_1)}{1+\lambda_1^1 \Delta g(t_1)}+
\alpha_1^2\,\int_{(t_1,t)} \frac{ \exp_g(\lambda_2^1;\cdot)}{\left(1+\lambda_2^2\, \Delta g\right)\,\exp_g(w_0^2\,\Delta g;\cdot)}\, f\, \operatorname{d}\mu_g\\
&+\alpha_2^2\,
\int_{(t_1,t)}
 \frac{\exp_g(\lambda_2^2;\cdot)}{\left(1+\lambda_2^1\, \Delta g\right)\,\exp_g(w_0^2\,\Delta g;\cdot)}\, f\, \operatorname{d}\mu_g
\Bigg]\\
&-\frac{\left(\alpha_1^2\, \exp_g(\lambda_2^1;t)
+\alpha_2^2\,\exp_g(\lambda_2^2;t) \right)}{\lambda_1^1-\lambda_1^2}
\Bigg[
\int_{[0,t_1)} \exp_g(\lambda_1^2;\cdot)^{-1} 
\frac{f}{1+\lambda_1^2 \Delta g}
\, \operatorname{d} \mu_g\\
&+\exp_g(\lambda_1^2;t_1)^{-1} 
\frac{f(t_1)\, \Delta g(t_1)}{1+\lambda_1^2 \Delta g(t_1)}
+\alpha_1^1\,\int_{(t_1,t)}\frac{ \exp_g(\lambda_2^1;\cdot)}{\left(1+\lambda_2^2\, \Delta g\right)\,\exp_g(w_0^2\,\Delta g;\cdot)}\, f\, \operatorname{d}\mu_g\\
&+
\alpha_2^1\,\int_{(t_1,t)} \frac{\exp_g(\lambda_2^2;\cdot)}{\left(1+\lambda_2^1\, \Delta g\right)\,\exp_g(w_0^2\,\Delta g;\cdot)}\, f\, \operatorname{d}\mu_g\Bigg],
\end{aligned}
\end{displaymath}
which finishes the proof.
\end{proof}

In the following example we obtain an explicit solution of 
problem~\eqref{eq:piecewisehA}-\eqref{eq:piecewisehB} for a particular 
choice of a derivator $g$.

\begin{exa} Let be $x_0,\,v_0\in \mathbb{C}$, $w_1,\,w_2\in \mathbb{R}$, $T>1$ and consider the derivator,
	\begin{displaymath}
		g(t)=\left\{
		\begin{array}{ll}
			t,\quad & t\le 1,\\
			t+\delta,\quad & t>1,
		\end{array}
		\right.
	\end{displaymath}
for some $\delta>0$. It is clear that $D_g=\{1\}$ and $\Delta g(1)=
\delta$. Moreover, given an element $\lambda \in \mathbb{C}$ 
such that $1+\lambda\,\delta \neq 0$, we have that 
\begin{displaymath}
	\exp_g(\lambda;t)=\left\{
	\begin{array}{ll}
		e^{\lambda\, t},\quad & 0\le t\le 1,\\
		\left(1+\lambda\,\delta\right) \, e^{\lambda\, t},\quad & 1<t\le T,
	\end{array}
	\right.
\end{displaymath}
In this case, the solution of system~\eqref{eq:systemalphak} is
\begin{displaymath}
\begin{aligned}
\alpha_1^1=&\frac{ e^{i\,(w_{1}-w_{2})}\,\left(w_{1}+w_{2}\right)\,\left(1+i\,\delta\,w_{1}\right)}{2\,w_{2}\,\left(1+i\,\delta\,w_{2}\right)},&
\alpha_2^1&=\frac{ e^{i\,(w_1+w_2) }\,\left(w_{2}-w_{1}\right)\,\left(1+i\,\delta\,w_{1}\right)}{2\,w_{2}\,\left(1-
i\,\delta\,w_{2}\right)},\\
\alpha_1^2=&\frac{ e^{-i\,(w_1+w_2) }\,\left(w_{2}-w_{1}\right)\,\left(1-i\,\delta\,w_{1}\right)}{2\,w_{2}\,\left(1+i\,\delta\,w_{2}\right)},&
\alpha_2^2&=\frac{ e^{i\,(w_2-w_1)}\,\left(w_{1}+w_{2}\right)\,\left(1-i\,\delta\,w_{1}\right)}{2\,w_{2}\,\left(1-i\,\delta\,w_{2}\right)}.
\end{aligned}
\end{displaymath}
Therefore, the solution of the corresponding 
problem~\eqref{eq:piecewisehA}-\eqref{eq:piecewisehB} is given 
by the following expression:
\begin{displaymath}
\begin{aligned}
v_h(t)=&\left(\frac{w_1\, x_0-i\,v_0}{2\,w_1} \right)\,
\bigg[ e^{i\,w_1\,t}\, \chi_{[0,1]}(t)+\bigg(
\frac{ e^{i\,(w_{1}-w_{2})}\,\left(w_{1}+w_{2}\right)\,\left(1+i\,\delta\,w_{1}\right)}{2\,w_{2}\,
}\, 
e^{i\,w_2\,t}\\
&
+\frac{ e^{i\,(w_1+w_2) }\,\left(w_{2}-w_{1}\right)\,\left(1+i\,\delta\,w_{1}\right)}{2\,w_{2}}\,e^{-i\,w_2\,t}\bigg)\,\chi_{(1,T]}(t)\bigg]\\
&+\bigg(\frac{w_1\,x_0+i\,v_0}{2\,w_1} \bigg)\, 
\bigg[
e^{-i\,w_1\,t}\, \chi_{[0,1]}(t)
+\bigg(
\frac{ e^{-i\,(w_1+w_2) }\,\left(w_{2}-w_{1}\right)\,\left(1-i\,\delta\,w_{1}\right)}{2\,w_{2}}
\,e^{i\,w_2\,t}\\
&
+\frac{ e^{-i\,(w_1-w_2)}\,\left(w_{1}+w_{2}\right)\,\left(1-i\,\delta\,w_{1}\right)}{2\,w_{2}}\,e^{-i\,w_2\,t} \bigg)\,\chi_{(1,T]}(t)\bigg]\\
=&\bigg[\left(\frac{w_1\, x_0-i\,v_0}{2\,w_1} \right)\,e^{i\,w_1\,t}+
 \left( \frac{w_1\,x_0+i\,v_0}{2\,w_1}\right)\,
 e^{-i\,w_1\,t}
 \bigg]\, \chi_{[0,1]}(t)\\
 &+\bigg[
 \left(\frac{w_1\, x_0-i\,v_0}{2\,w_1} \right)\,
 \frac{ e^{i\,(w_{1}-w_{2})}\,\left(w_{1}+w_{2}\right)\,\left(1+i\,\delta\,w_{1}\right)}{2\,w_{2}\,
}\, 
e^{i\,w_2\,t}\\
&+\left( \frac{w_1\,x_0+i\,v_0}{2\,w_1}\right)\,
\frac{ e^{-i\,(w_1-w_2)}\,\left(w_{1}+w_{2}\right)\,\left(1-i\,\delta\,w_{1}\right)}{2\,w_{2}}\,e^{-i\,w_2\,t}
\bigg]\, \chi_{(1,T]}(t)\\
&+\bigg[
\left(\frac{w_1\, x_0-i\,v_0}{2\,w_1} \right)\,
\frac{ e^{i\,(w_1+w_2) }\,\left(w_{2}-w_{1}\right)\,\left(1+i\,\delta\,w_{1}\right)}{2\,w_{2}}\,e^{-i\,w_2\,t}\\
&+
\left( \frac{w_1\,x_0+i\,v_0}{2\,w_1}\right)\,
\frac{ e^{-i\,(w_1+w_2) }\,\left(w_{2}-w_{1}\right)\,\left(1-i\,\delta\,w_{1}\right)}{2\,w_{2}}
\,e^{i\,w_2\,t}
\bigg]\, \chi_{(1,T]}(t).
\end{aligned}
\end{displaymath}
Observe that, given $t\in [0,T]$, 
\begin{displaymath}
\begin{aligned}
\widetilde{v}_h(t):=\lim_{\delta \to 0} v_h(t)=&\bigg[\left(\frac{w_1\, x_0-i\,v_0}{2\,w_1} \right)\,e^{i\,w_1\,t}+
 \left( \frac{w_1\,x_0+i\,v_0}{2\,w_1}\right)\,
 e^{-i\,w_1\,t}
 \bigg]\, \chi_{[0,1](t)}\\
 &+\bigg[
 \left(\frac{w_1\, x_0-i\,v_0}{2\,w_1} \right)\,
 \frac{ e^{i\,(w_{1}-w_{2})}\,\left(w_{1}+w_{2}\right)}{2\,w_{2}\,
}\, 
e^{i\,w_2\,t}\\
&+\left( \frac{w_1\,x_0+i\,v_0}{2\,w_1}\right)\,
\frac{ e^{-i\,(w_1-w_2)}\,\left(w_{1}+w_{2}\right)}{2\,w_{2}}\,e^{-i\,w_2\,t}
\bigg]\, \chi_{(1,T]}(t)\\
&+\bigg[
\left(\frac{w_1\, x_0-i\,v_0}{2\,w_1} \right)\,
\frac{ e^{i\,(w_1+w_2) }\,\left(w_{2}-w_{1}\right)}{2\,w_{2}}\,e^{-i\,w_2\,t}\\
&+
\left( \frac{w_1\,x_0+i\,v_0}{2\,w_1}\right)\,
\frac{ e^{-i\,(w_1+w_2) }\,\left(w_{2}-w_{1}\right)}{2\,w_{2}}
\,e^{i\,w_2\,t}
\bigg]\, \chi_{(1,T]}(t)
\end{aligned}
\end{displaymath}
is such that $\widetilde{v}_h\in\{v\in\mathcal{C}^1([0,T];\mathbb{R}): v'\in \mathcal{AC}([0,T];\mathbb{R})\}$ is the solution of the 
classical one-dimensional linear Helmholtz equation with 
piecewise-constant coefficients
\begin{displaymath}
(\widetilde{v}_h)'(t)=v_0+\int_{[0,t]} -w_0^2\, 
\widetilde{v}_h\, \operatorname{d} m,\quad t\in [0,T],
\end{displaymath}
where $m$ is the Lebesgue measure in $[0,T]$ and $\widetilde{v}_h(0)=x_0$. Observe that, for $\delta=0$, the 
$g$-derivative coincides with the usual derivative. It is also remarkable 
that $\widetilde{v}_h$ and $(\widetilde{v}_h)'$ are continuous 
functions, whereas the second derivative $(\widetilde{v}_h)''$ 
undergoes jumps at the points of discontinuity of $w_0$ 
(see \cite{BARUCH2007820}). This behavior can be seen in 
Figures~\ref{fig1}, \ref{fig2} and \ref{fig3}.

\begin{figure}[H]
\centering
\includegraphics[width=0.7\linewidth]{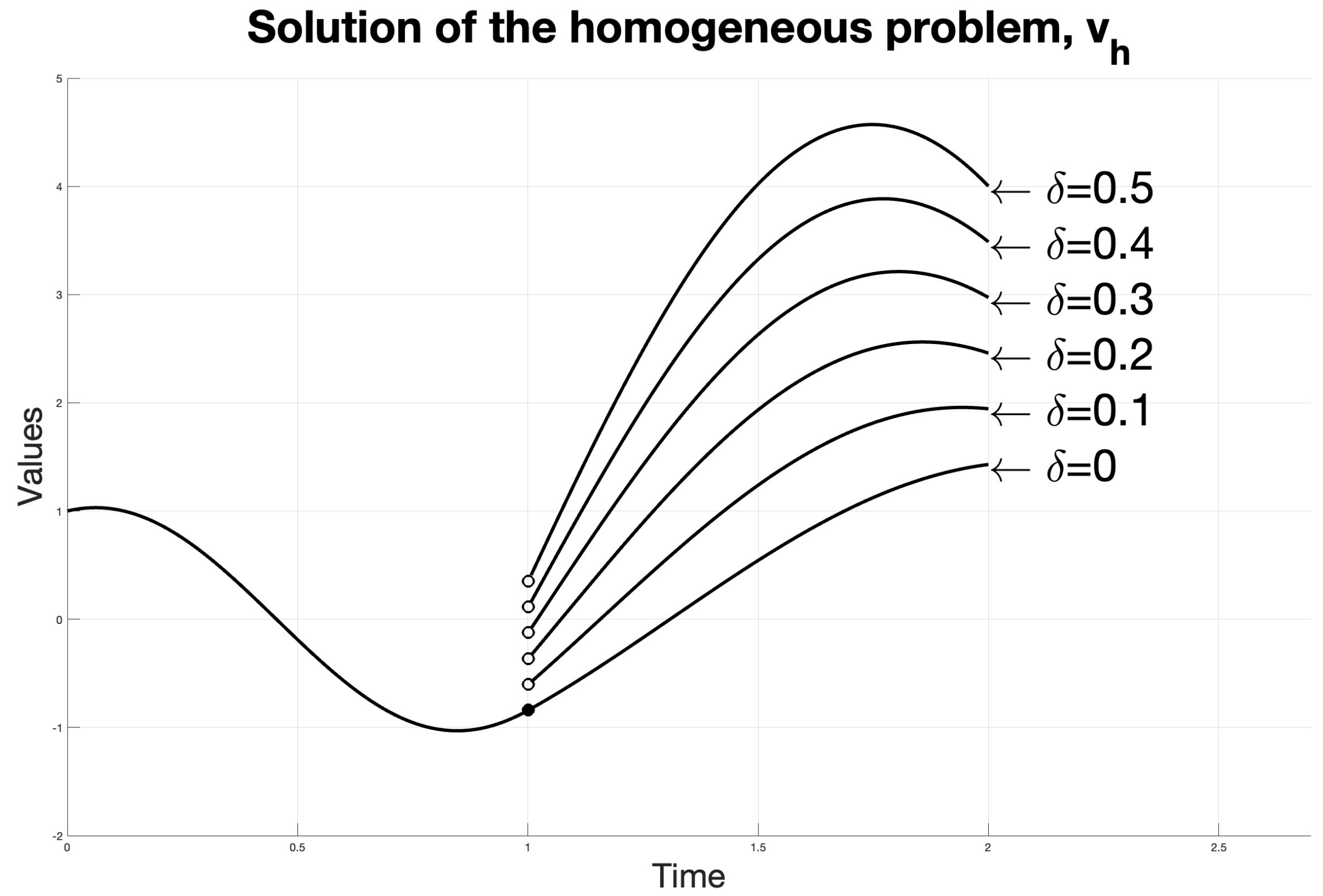}
\caption{Solution of the homogeneous problem~\eqref{eq:piecewisehA}-\eqref{eq:piecewisehB} for different values of $\delta$. Observe that for $\delta=0$ the $g$-continuity is equivalent to the continuity in the usual sense.}
\label{fig1}
\end{figure}

\begin{figure}[H]
\centering
\includegraphics[width=0.7\linewidth]{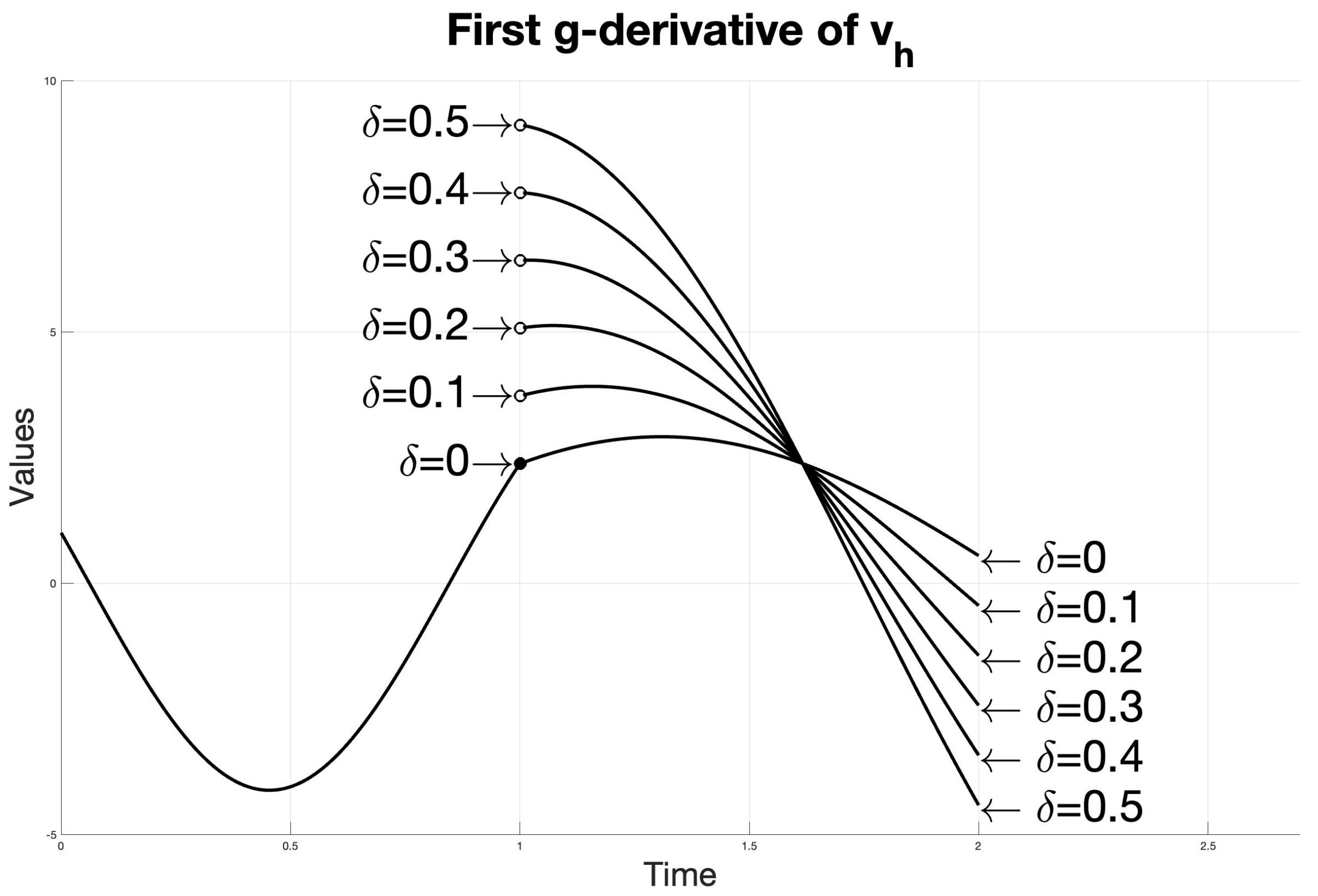}
\caption{First $g$-derivative of the solution of~\eqref{eq:piecewisehA}-\eqref{eq:piecewisehB} for different values of $\delta$. Observe that for $\delta=0$ the $g$ derivative is equivalent to the usual derivative and the first derivative of 
the solution is continuous in the usual sense.}
\label{fig2}
\end{figure}

\begin{figure}[H]
\centering
\includegraphics[width=0.7\linewidth]{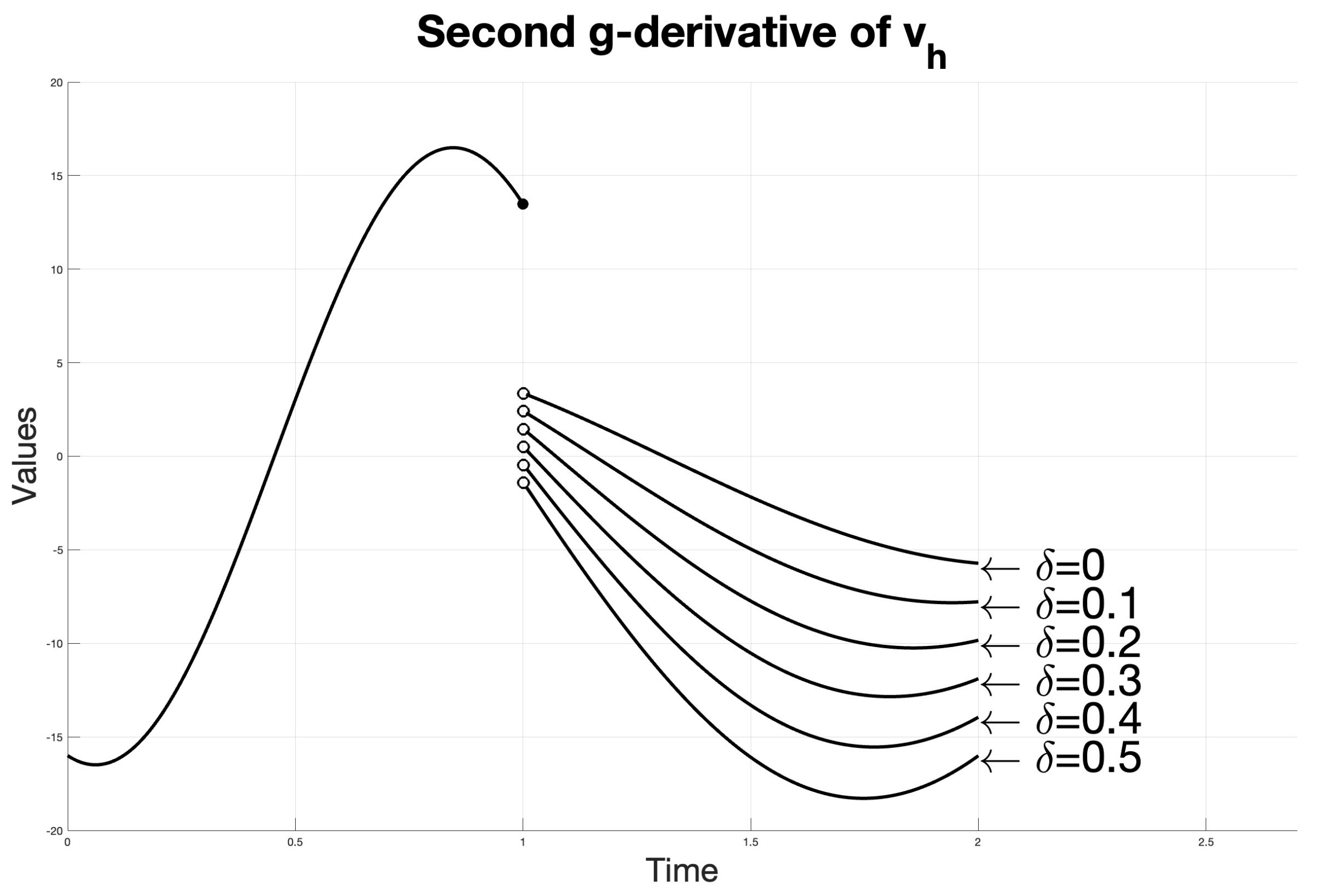}
\caption{Second $g$-derivative of the solution of~\eqref{eq:piecewisehA}-\eqref{eq:piecewisehB} for different values of $\delta$. Observe that the second $g$-derivative is $g$-continuous 
for all $\delta>0$, whereas for $\delta=0$ the second 
derivative $(\widetilde{v}_h)''=(\widetilde{v}_h)''_g$ has a discontinuity point at $t=1$.}
\label{fig3}
\end{figure}
\end{exa}

\begin{rem} This last example suggests that there is some form of continuity of the solution of problem~\eqref{eq:piecewisehA}-\eqref{eq:piecewisehB} with respect to $g$. Indeed, this is in general the case. Let $\mathcal{B}$ be the Borel $\sigma$-algebra associated to the usual topology of $[a,b]\subset\mathbb R$ and let $\mathcal{M}([a,b], \mathcal{B})$ be the space of all signed measures of bounded variation defined on $\mathcal{B}$. The total variation $\|\cdot\|$  is a norm on $\mathcal{M}([a,b], \mathcal{B})$ with which this space is a Banach space \cite[Theorem 4.6.1]{Bogachev}.  If we now fix a bounded Borel function $f:[0,T]\to\mathbb{C}$ then the map $T:(\mathcal{M}([a,b], \mathcal{B}),|\cdot|)\to \mathbb{C}$ such that $T\mu=\int f\, \operatorname{d} \mu$ is a bounded linear functional. Indeed, for every $\mu\in \mathcal{M}([a,b], \mathcal{B})$,
\[|T\mu|=\left|\int f\operatorname{d} \mu\right|\le\int |f|\operatorname{d} |\mu|\le \|f\|_\infty\|\mu\|.\]
\end{rem}

\section*{Acknowledgments}
The authors would like to thank the anonymous referee for their comments, suggestions and corrections, as they have greatly contributed to improve the quality of the manuscript.

The authors were partially supported by Xunta de Galicia, project ED431C 2019/02, and by the Agencia Estatal de Investigaci\'on (AEI) of Spain under grant MTM2016-75140-P, co-financed by the European Community fund FEDER.

\bibliography{refs-gso2}

\end{document}